\documentclass[reqno,10pt]{amsart}

\usepackage{a4wide}
\usepackage{color}
\usepackage[T2A]{fontenc}
\usepackage[utf8]{inputenc}
\usepackage[english]{babel}
\numberwithin{equation}{section}
\usepackage[colorlinks,citecolor=blue,linkcolor=red]{hyperref}
\usepackage[noabbrev,capitalise,nameinlink]{cleveref}
\hypersetup{colorlinks={true},linkcolor={red},citecolor=blue}
\usepackage{csquotes}

\makeatletter
\newcommand{\bg}{\bBigg@{0.5}}
\newcommand{\biggg}{\bBigg@{3}}
\newcommand{\Biggg}{\bBigg@{3.5}}
\newcommand{\bigggg}{\bBigg@{4}}
\newcommand{\bdot}{\text{$\,\cdot\,$}}
\newcommand{\lims}{\varlimsup}
\DeclareMathOperator*{\esup}{ess\,sup}
\DeclareMathOperator*{\einf}{ess\,inf}
\newcommand{\limi}{\varliminf}
\newcommand{\supp}{\operatorname{supp}}

\newcommand{\card}{\operatorname{card}}
\newcommand{\D}{\mathrm{d}}

\newcommand{\E}{\mathrm{E}}
\newcommand{\Leb}{\mathrm{L}}
\newcommand{\BLIP}{\mathrm{BLIP}}

\newcommand{\dist}{\operatorname{dist}}
\newcommand{\R}{\mathbb{R}}

\newcommand{\Me}{\mathfrak{M}}
\newcommand{\X}{\mathsf{X}}
\newcommand{\Y}{\mathsf{Y}}
\newcommand{\dY}{\mathsf{d}_{\mathsf{Y}}}
\newcommand{\Z}{\mathsf{Z}}
\newcommand{\dZ}{\mathsf{d}_{\mathsf{Z}}}
\newcommand{\B}{\mathrm{B}}
\newcommand{\A}{\mathrm{A}}
\newcommand{\m}{\mathfrak{m}}

\newcommand{\dX}{\mathsf{d}}

\newcommand{\Lip}{\mathrm{Lip}}
\newcommand{\LIP}{\mathrm{LIP}}
\newcommand{\lip}{\mathrm{lip}}
\newcommand{\Q}{\mathsf{Q}}
\renewcommand{\P}{\mathcal{P}}

\newcommand*\rel@kern[1]{\kern#1\dimexpr\macc@kerna}
\newcommand*\widebar[1]{%
  \begingroup
  \def\mathaccent##1##2{%
    \rel@kern{0.8}%
    \overline{\rel@kern{-0.8}\macc@nucleus\rel@kern{0.2}}%
    \rel@kern{-0.2}%
  }%
  \macc@depth\@ne
  \let\math@bgroup\@empty \let\math@egroup\macc@set@skewchar
  \mathsurround\z@ \frozen@everymath{\mathgroup\macc@group\relax}%
  \macc@set@skewchar\relax
  \let\mathaccentV\macc@nested@a
  \macc@nested@a\relax111{#1}%
  \endgroup
}

\newcommand{\wbar}{\widebar}

\makeatother

\usepackage{amsmath, amsfonts, amssymb, amsthm, mathtools, bbm, esint, mathrsfs, nccmath}
\usepackage[backend=biber,style=numeric-comp,sortcites,giveninits,doi=false,isbn=false,url=false,eprint=false, maxnames=10]{biblatex}

\usepackage{resmes}

\AtBeginBibliography{\footnotesize\sloppy}

\DeclareNameAlias{author}{family-given}

\theoremstyle{plain}
\newtheorem{theo}{Theorem}[section]
\newtheorem{lem}[theo]{Lemma}
\newtheorem{prop}[theo]{Proposition}
\newtheorem{cor}[theo]{Corollary}

\theoremstyle{definition}
\newtheorem{defi}[theo]{Definition}

\theoremstyle{remark}
\newenvironment{rema} {\begin{proof}[\rm{\textbf{Remark}}]} {\end{proof}}

\defbibheading{bibliography}[\refname]{%
  \section*{#1}
}

\usepackage{relsize}

\title[On estimates of the BBM type in the singular context]{On estimates\\
of the Bourgain--Brezis--Mironescu type\\
in the singular context}

\author[R. D. Oleinik]{Roman D. Oleinik}

\address{{International School for Advanced Studies, Trieste, Italy}\newline
\indent {Moscow Institute of Physics and
Technology, Dolgoprudny, Russia}\newline
\indent {Steklov Mathematical Institute of the Russian Academy of Sciences, Moscow, Russia}
}

\email{oleinik.r@phystech.edu}

\begin{document}
\allowdisplaybreaks

\begin{abstract}
We extend the characterization of the Bourgain--Brezis--Mironescu type for maps from certain metric measure spaces to arbitrary metric spaces to a broader class of mollifiers.
\end{abstract}

\maketitle
\tableofcontents
\newpage
\section{Introduction}

\subsection{Overview} Within the seminal paper \cite{BBM01} by J. Bourgain, H. Brezis, and P. Mironescu, there was established one remarkable property of Sobolev functions on Euclidean domains expressed in terms of the limiting behavior of a special class of nonlocal functionals. Namely, the following was obtained. Let $O\subseteq \mathbb{R}^d$, $d\in \mathbb{N}$, be a bounded, smooth domain and let $(\varrho_{\delta})_{\delta\in (0,1)}$ be a family of radial mollifiers on $\mathbb{R}^d$, understood as $\mathcal{L}^1$-measurable functions $[0,+\infty)\to [0,+\infty]$, where $\mathcal{L}^1$ denotes the Lebesgue measure on $\mathbb{R}$, such that
\begin{equation}\label{eq:OrigNormalCond}
    \int\limits_{\bg[0,+\infty\bg)} \varrho_{\delta}(t) t^{d-1}\D \mathcal{L}^1(t)=1
\end{equation}
for every $\delta\in (0,1)$ and
\begin{equation}\label{eq:OrigLocalCond}
    \lims\limits_{\delta \searrow 0}\int\limits_{\bg(r,+\infty\bg)} \varrho_{\delta}(t) t^{d-1}\D \mathcal{L}^1(t)=0
\end{equation}
for every $r\in (0,+\infty)$. Then the formula below holds for any $u\in \Leb^p(O)$ with $p\in (1,+\infty)$:
\begin{equation}
\label{eq:OrigRes}
\begin{gathered}
\lim\limits_{\delta\searrow 0}\int\limits_{O\times O} \Bigg(\frac{\big|u(x')-u(x)\big|}{\big\|x'-x\big\|_d}\Bigg)^p\varrho_{\delta}\Big(\big\|x'-x\big\|_d\Big) \D \big(\mathcal{L}^d \otimes\mathcal{L}^d\big)(x,x')=\\
=C(p,d)\mathrm{E}_{p}[u](O),
\end{gathered}
\end{equation}
where $\|\bdot\|_d$ and $\mathcal{L}^d$ are the Euclidean norm and the Lebesgue measure, respectively, on $\mathbb{R}^d$, $C(p,d)$ is a positive constant depending only on $p$ and $d$, and $\mathrm{E}_{p}[u](O)$ is the $p$-energy of $u$ on $O$ given by
\begin{equation}
\mathrm{E}_{p}[u](O)\coloneqq \begin{dcases}\int\limits_{O} \Big(\big\|(\nabla u)(x)\big\|_d\Big)^p \D \mathcal{L}^d(x),\quad u\in \mathrm{W}^{1,1}_{Loc}(O),\\
+\infty, \quad \text{otherwise},
\end{dcases}
\end{equation}
where $\nabla u$, in turn, denotes the weak gradient of $u$. As it was soon discovered, an analogous feature is inherent in functions of bounded variation as well. Namely, it was shown by J. D{\'a}vila in \cite{D02} that \eqref{eq:OrigRes}, with $p=1$, holds for any $u\in \Leb^1(O)$, if one understands $\mathrm{E}_{1}[u](O)$, the $1$-energy of $u$ on $O$, as the total variation of $u$ on $O$. These results led eventually to the appearance in the literature of various similar statements, which are now typically called BBM-type ones after
Bourgain--Brezis--Mironescu.

In last decades, one of the major trends in analysis lies in the extension of the classical first-order calculus to the singular context of so-called metric measure spaces, usually understood as a separable metric space $(\X,\dX)$ endowed with a locally finite, Borel outer measure $\m$. As a part of this process, there began to appear numerous researches presenting BBM-type results in such a setting. And a prototypical problem arising therein, more natural for the less regular framework, can be formulated as follows: to find suitable sufficient conditions on the mollifiers under which one can achieve an equivalence, up to some universal structural constants, between the right-hand side and the lower and upper limits as $\delta\searrow 0$ of the integral in the left-hand side of \eqref{eq:OrigRes}, with all the objects appropriately generalized to the new situation. Statements of such a form are occasionally called the BBM-type characterization. The very first contribution in this direction was made by V. Munnier in \cite{M15}, where the author dealt with a very special family of mollifiers, namely the one related to the fractional Sobolev seminorms. With an analogous family involved, the full-fledged solution to the problem was lately obtained by S. Di Marino and M. Squassina in \cite{DS19}. Within the subsequent paper \cite{LPZ22} by P. Lahti, A. Pinamonti, and X. Zhou, the class of admissible families was significantly enlarged so that it included quite many of mollifiers of interest. And essentially the same families fell into the scope of the later work \cite{O25} by the author of the current manuscript, where, among other things, all the above-described subject was moreover adapted to maps with values in arbitrary metric spaces.

At the same time, the conditions on mollifiers considered in \cite{LPZ22} and in \cite{O25} did not incorporate all relevant kinds of such. The main restriction of those, apart of other less important ones, can be loosely summarized as follows: the functions
\begin{equation}
    (0,+\infty)\to [0,+\infty]\colon t\mapsto \frac{\varrho_{\delta}(t)}{t^p},\quad \delta\in (0,1),
\end{equation}
must be non-increasing. A meaningful example of mollifiers not fitting into this pattern is the family with the functions given by
\begin{equation}\label{eq:NonMonFamEx}
    \varrho_{\delta}(t)\coloneqq \begin{dcases}
        \frac{1}{\big|\ln(\delta)\big| t^d}, \quad t\in (\delta,1],\\
        0, \quad \text{otherwise},
    \end{dcases}
\end{equation}
for every $\delta\in (0,1)$. Although these functions violate the above-mentioned condition ``at just one point'', the techniques exploited in the previous articles were not adequate to encompass this family. In correspondence with this, our goal for this paper is to eliminate the presented drawback, namely we want to provide a version of the BBM-type characterization that is applicable to even those mollifiers that are not necessarily ``monotone'' in the above sense, including the specific family given by \eqref{eq:NonMonFamEx}. And in line with the approach developed in \cite{O25}, we are to deal with metric-valued maps, not just with real-valued functions, within our work.

\subsection{Results} The main subject of our study can be formulated as follows. Given a metric measure space $(\X, \mathsf{d}, \mathfrak{m})$, an exponent $p\in [1,+\infty)$, a family $(\rho_{\delta})_{\delta\in (0,1)}$ of functions $\X\times \X\to [0,+\infty]$, a metric space $(\Y, \dY)$, an open set $O\subseteq \X$, and a map $f\colon O\to \Y$, we will analyze the limiting behavior as $\delta\searrow 0$ of the quantities
\begin{gather}
\label{eq:MainQuant}
\int\limits_{O\times O} \Bigg(\frac{\dY\big(f(x),f(x')\big)}{\dX(x,x')}\Bigg)^p\rho_{\delta}(x,x') \D (\mathfrak{m}\otimes \mathfrak{m})(x,x') , \quad\text{$\delta\in (0,1)$},
\end{gather}
under certain assumptions on all the objects involved, which are to be indicated later. More specifically, we will seek for mutual quantitative estimates between the lower and upper limits as $\delta \searrow 0$ of the quantities in \eqref{eq:MainQuant} and $\E_p[f](O)$, the $p$-energy of $f$ on $O$, which is understood in the sense of \cref{def:MetrCh}.

Within the paper, the source metric measure space will always be meant locally complete. Moreover, in line with all the previous works on the topic, we will use certain variations of the so-called doubling and Poincar{\'e} conditions as the main structural assumptions on the underlying space. And we incorporate them in the concept of a $p$-Poincar{\'e} space, the meaning of which is disclosed in \cref{def:PoinSpace}. In turn, the requirements for mollifiers that we will deal with are listed within \cref{def:AdmFam}, and we will refer to families of mollifiers fulfilling these as $p$-admissible. At last, we will need to impose on $O$ a certain extension property, which we formulate in \cref{def:StrExtDom}. 

Our main result, in its most digestible form, can be stated as follows.
\begin{theo}\label{theo:Res}
Let $(\X, \dX,\m)$ be a locally complete metric measure space. Let $p\in [1,+\infty)$. Suppose $(\X, \dX,\m)$ is a $p$-Poincar{\'e} space. Let $(\rho_{\delta})_{\delta\in (0,1)}$ be a $p$-admissible family. Then there exists a universal constant $C\in (0,+\infty)$ such that the following holds. Let $(\Y,\dY)$ be a metric space, let $O\subseteq \X$ be an open set having the strong $p$-extension property with respect to $\Y$, let $f\in \Leb^p(O, \Y)$. Then the estimates below take place:
\begin{equation}
\label{eq:ResIneq}
    \begin{gathered}
\frac{1}{C}\mathrm{E}_p[f](O)\leq\\
\leq \limi\limits_{\delta\searrow 0} \int\limits_{O\times O} \Bigg(\frac{\dY\big(f(x),f(x')\big)}{\dX(x,x')}\Bigg)^p\rho_{\delta}(x,x') \D (\mathfrak{m}\otimes \mathfrak{m})(x,x') 
\leq  \\
\leq\lims\limits_{\delta\searrow 0} \int\limits_{O\times O} \Bigg(\frac{\dY\big(f(x),f(x')\big)}{\dX(x,x')}\Bigg)^p\rho_{\delta}(x,x') \D (\mathfrak{m}\otimes \mathfrak{m})(x,x')  \leq \\
\leq C\mathrm{E}_p[f](O).
\end{gathered}
\end{equation}
\end{theo}
\noindent This statement repeats \cite[Theorem 1]{O25} essentially word for word, with the only relevant difference being in that here the class of admissible mollifiers is understood in a broader sense in comparison with the earlier result. And whereas we do not justify this latter claim in our paper explicitly, we insist it to be true. As an illustration for this, it will be demonstrated in \cref{ss:ExFamMol} that the family given by \eqref{eq:NonMonFamEx} is indeed covered by \cref{theo:Res}, which was exactly our initial motivation.

The proof of the theorem will be carried out in two stages, separately for the lower and upper bounds. The corresponding parts can be found in \cref{theo:LowBoundTh} and \cref{theo:UpBoundTh}, respectively. We also draw the reader's attention to that both of them will be formulated in slightly more general terms, adapted specifically to each of the statements. And to reproduce the main result one just needs to combine the given theorems with \cref{prop:DoubStrDoub}, \cref{def:PoinSpace}, and \cref{def:AdmFam}.

The ideas we follow to prove our result are to some extent similar to those that we exploited within \cite{O25}. More specifically, the proof of the upper inequality in \eqref{eq:ResIneq} uses again a pointwise inequality implied by the doubling and Poincar{\'e} conditions, but stated in a slightly stronger form. In turn, the proof of the lower inequality in \eqref{eq:ResIneq} relies on an approximation procedure, still based on discrete convolutions, that is somewhat unusual in comparison with the more standard approach appearing in the literature.

The organization of the paper is as follows. We formulate all the necessary definitions in \cref{ss:Prelim}, including \cref{ss:GenerNot} with a general notation, \cref{ss:SobolCal} with a relevant information on Sobolev calculus, \cref{ss:DoubPonCon} with notions related to doubling and Poincar{\'e} conditions, and \cref{ss:ReqMol} with requirements for mollifiers. And all the necessary statements are given within \cref{ss:Proofs}, consisting of \cref{ss:LowBound} with the proof of the lower-bound inequality, of \cref{ss:UpBound} with the proof of the upper-bound inequality, of \cref{ss:SimplAssMol} with simplified versions of the requirements, and of \cref{ss:ExFamMol} with several specific examples of mollifiers.

\section{Preliminaries}\label{ss:Prelim}

\subsection{General notation}\label{ss:GenerNot}
We start by introducing some notation.

Until the end of the manuscript, the following data is assumed fixed:
\begin{itemize}
\item[$\bullet$] a locally complete, separable metric space $(\X, \dX)$;
\item[$\bullet$] a locally finite, Borel outer measure $\m$ on $(\X, \dX)$.
\end{itemize}
The resulting metric measure space $(\X,\dX, \m)$ will be denoted just by $\X$.

To avoid discussing separately degenerate cases, we initially assume that our space $\X$ is of positive diameter and of positive measure.

By $\supp(\m)$ we will denote the support of $\m$, defined as the set of all $x\in \X$ such that any neighborhood of $x$ has positive measure. It is known that $\supp(\m)$ is a closed subset of $\X$ and that $\m\big(\X\backslash \supp(\m)\big)=0$. And provided that $\supp(\m)=\X$, we will say that $\X$ has full support.

We use the notation $\mathbb{R}_{\geq \alpha}\coloneqq [\alpha,+\infty)$ and $\mathbb{R}_{>\alpha}\coloneqq (\alpha,+\infty)$ for any given $\alpha\in [0,+\infty)$. We will also write $\overline{\mathbb{R}}_{\geq 0}$ instead of $[0,+\infty]$.

Given an at most countable set $S$, by $\card(S)$ we mean the cardinality of $S$, understood as an element of $\mathbb{N}_0\cup\{+\infty\}$.

In what follows, we adopt the standard conventions that $0\cdot (+\infty)=\tfrac{0}{0}=0$.

Given $S_1,S_2\subseteq \X$, let $\dist(S_1,S_2)$ be the distance between $S_1$ and $S_2$ with respect to $\dX$. Here we use the convention that the distance between any subset of $\X$ and $\varnothing$ is infinite.

Given $x\in \X$, $r\in \mathbb{R}_{> 0}$, $r'\in (0,r)$, $\tau\subseteq \R_{\geq 0}$, and $S\subseteq \X$, we introduce the following sets:
\begin{gather}
\mathrm{B}(x,r)\coloneqq \big\{x'\in \X \mid \mathsf{d}(x,x')\leq r\big\}, \\
\mathrm{B}(S,r)\coloneqq \Big\{x'\in \X \bigm| \dist\big(S,\{x'\}\big)\leq r\Big\},\\
\A(x,\tau)\coloneqq \big\{x'\in \X\mid \dX(x,x')\in \tau\big\},\\
\mathrm{A}(x,r,r')\coloneqq \A\big(x,(r',r]\big),\\ \mathrm{A}(x,\infty,r)\coloneqq \A\big(x,(r,+\infty)\big).
\end{gather}

Given $S_0,S\subseteq \X$ with $S\subseteq S_0$, we write $S\Subset S_0$ whenever $\dist\big(S,\X\backslash S_0\big)>0$.

Given $S\subseteq \X$, by $\chi_{S}$ we mean the characteristic function of $S$. Given arbitrary sets  $S_1,S_2,S$ with $S\subseteq S_1$ and a map $f\colon S_1\to S_2$, by $f\big|_{S}$ we denote the restriction of $f$ to $S$.

Given a metric space $(\Y,\dY)$ and a set $S$, a map $f\colon S\to \Y$ is called separably valued if its image is separable as a subspace of $\Y$.

Let $p\in \R_{\geq 1}$, let $(\Y,\dY)$ be a metric space, let $E\subseteq \X$ be a $\m$-measurable set, let $S\subseteq \Y$. By $\Me(E,S)$, we mean the collection of all $\m$-measurable maps $f\colon E\to S$ such that, for some $E'\subseteq E$ with $\m\big(E\backslash E'\big)=0$, the map $f\big|_{E'}$ is separably valued. Then, by $\Leb^p(E)$ we mean the collection of all functions belonging to $\Me(E,\mathbb{R})$ that are $p$-integrable on $E$. Next, let $\Leb^p_{Loc}(E)$ denote the collection of all functions $u\in \Me(E, \mathbb{R})$ with the property that, given any $x\in E$, there is $r\in \mathbb{R}_{>0}$ such that $u\in \Leb^p\big(E\cap \mathrm{B}(x,r)\big)$. Finally, let $\Leb^p(E,\Y)$ stand for the collection of all maps $f\in \Me(E,\Y)$ for which there is $y\in \Y$ such that $\mathsf{d}_{\Y}\big(f(\bdot),y\big)\in \Leb^p(E)$.

We will also write $\Me_+\big(\X\times \X\big)$ for the collection of all $(\m\otimes \m)$-measurable functions $\X\times \X\to \overline{\mathbb{R}}_{\geq 0}$.

Given a $\m$-measurable set $E\subseteq \X$ and $u\in \Leb^1(E)$, we introduce the quantity
\begin{equation}
\langle u\rangle_E\coloneqq \fint\limits_{E} u(x')\D\mathfrak{m}(x')\coloneqq \begin{dcases}
     \tfrac{1}{\mathfrak{m}(E)}\int\limits_{E} u(x')\D\mathfrak{m}(x'),\quad \m(E)\in \R_{>0},\\
     0,\quad \text{otherwise}.
\end{dcases}
\end{equation}

Let $p\in \R_{\geq 1}$, let $O\subseteq \X$ be an open set, let $u\in \Me(O,\mathbb{R})$. A sequence $(v_n)_{n\in \mathbb{N}}\subseteq \Me(O,\mathbb{R})$ is said to converge to $u$ in $\Leb^p_{Loc}(O)$ as $n\to +\infty$ provided that for any $x\in O$ there is $r\in \mathbb{R}_{>0}$ with $\mathrm{B}(x,r)\subseteq O$ such that
\begin{equation}
    \lims\limits_{n\to +\infty}\int\limits_{\mathrm{B}(x,r)} \Big|v_n(x')-u(x')\Big|^p\D \mathfrak{m}(x')=0.
\end{equation}

Given a metric space $(\Y,\dY)$, the collection of all bounded $1$-Lipschitz functions $\Y\to \mathbb{R}$ is denoted by $\BLIP_1(\Y)$.

From here on, given a metric space $(\Y,\dY)$, a set $S\subseteq \X$, and a map $f\colon S\to \Y$, by $\dX_f$ and $\Q_f$ we mean the functions $\dX_f, \Q_f\colon S\times S\to \mathbb{R}_{\geq 0}$ defined as
\begin{gather}
\mathsf{d}_f(x,x')\coloneqq \dY\big(f(x),f(x')\big),\\
\Q_f(x,x')\coloneqq \begin{dcases}
    \frac{\dX_f(x,x')}{\dX(x,x')},\quad x\neq x',\\
    0, \quad \text{otherwise}.
\end{dcases}
\end{gather}
This notation will also be used occasionally to substitute the corresponding full expressions. 

The following simple statement provides a measure-theoretic property necessary for all the future reasoning concerning metric-valued maps.
\begin{prop}\label{prop:DistFuncMeas}
Let $(\Y,\dY)$ be a metric space, let $E\subseteq \X$ be a $\m$-measurable set, let $f\in \Me(E,\Y)$. Then $\dX_f$ and $\Q_f$ are $(\mathfrak{m}\otimes \mathfrak{m})$-measurable.
\begin{proof}
We give an outline of the corresponding proof only. One needs to find $E'\subseteq E$ with $\m\big(E\backslash E'\big)=0$ such that the map $f\big|_{E'}$ is separably valued, take then the distance functions with respect to the points from a countable dense subset of $f(E')$, and consider finally truncations of those functions. This almost immediately gives the desired statement via some standard measure-theoretic arguments.
\end{proof}
\end{prop}

Given $S\subseteq \X$, the global Lipchitz constant of a function $u\colon S\to \R$ is given by
\begin{equation}
    \Lip[u]\coloneqq \sup\limits_{x,x'\in S }\Q_u(x,x').
\end{equation}

We record the following trivial proposition.
\begin{prop}\label{prop:GlLipConstProp}
    Let $S\subseteq \X$. Then the following assertions hold:
    \begin{enumerate}
        \item[$\rm i)$] for all functions $u,v\colon S\to \mathbb{R}$, one has
        \begin{equation}
            \Lip[uv]\leq \sup(u)\Lip[v]+\sup(v)\Lip[u];
        \end{equation}
        \item[$\rm ii)$] for any function $u\colon S\to \R_{>0}$, one has
        \begin{equation}
            \Lip\Big[\tfrac{1}{u}\Big]\leq \frac{\Lip[u]}{\big(\inf(u)\big)^2}.
        \end{equation}
    \end{enumerate}
    \begin{proof}
        The proof is straightforward.
    \end{proof}
\end{prop}

Let $O\subseteq \X$ be an open set. By $\mathrm{LIP}^{Loc}(O)$ we mean the collection of all functions $u\colon O\to \R$ with the property that for any $x\in O$ there is $r\in \mathbb{R}_{>0}$ with $\mathrm{B}(x,r)\subseteq O$ such that $u\big|_{\mathrm{B}(x,r)}$ is Lipschitz. Given a function $u\colon O\to \R$, we define a function $\mathrm{lip}[u]\colon O\to \overline{\mathbb{R}}_{\geq 0}$, known as the local Lipschitz constant of $u$, by
\begin{gather}
\mathrm{lip}[u](x)\coloneqq \begin{dcases}
    \lims\limits_{x'\rightarrow x}\Q_u(x,x'),\quad \text{$x$ is a limit point of $O$},\\
    0, \quad \text{otherwise}.
\end{dcases}
\end{gather}
It is not hard to see that $\mathrm{lip}[u]$ is Borel for every $u\in \mathrm{LIP}^{Loc}(O)$.

Given $p\in \R_{\geq 1}$ and an open set $O\subseteq \X$, a point $x\in O\cap \supp(\m)$ is called a $p$-Lebesgue point of $u\in \Leb^1_{Loc}(O)$ if
\begin{equation}
    \lims\limits_{r\searrow 0}\fint\limits_{\B(x,r)}\big(\dX_u(x,x')\big)^p \D \m(x')=0.
\end{equation}

Later, we will need to use the generalized version of the dominated convergence theorem, which we state as follows.
\begin{prop}\label{prop:GenLebTh}
Let $E\subseteq \X$ be a $\m$-measurable set, let $\Gamma\subseteq \mathbb{R}$ be a set having zero as a limit point, let $(u_{\gamma})_{\gamma\in \Gamma},(v_{\gamma})_{\gamma\in \Gamma}\subseteq \Me\big(E,\mathbb{R}_{\geq 0}\big)$. Suppose the following conditions are satisfied:
\begin{itemize}
\item[$\rm i) $] for each $\gamma\in \Gamma$, one has $u_{\gamma}(x)\leq v_{\gamma}(x)$ for $\mathfrak{m}$-a.e. $x\in E$;
\item[$\rm ii) $] it holds that $(u_{\gamma})_{\gamma\in \Gamma}$ converges $\mathfrak{m}$-a.e. on $E$ to some function $u\in \Me\big(E,\mathbb{R}_{\geq 0}\big)$ as $\gamma \rightarrow 0$;
\item[$\rm iii) $] it holds that $(v_{\gamma})_{\gamma\in \Gamma}$ converges in $\Leb^1(E)$ to some function $v\in \Leb^1(E)$ as $\gamma \rightarrow 0$.
\end{itemize}
Then it holds that $u\in\Leb^1(E)$ and that $(u_{\gamma})_{{\gamma}\in \Gamma}$ converges in $\Leb^1(E)$ to $u$ as $\gamma \rightarrow 0$.
\begin{proof}
We refer to \cite[Lemma 3.14]{GT21} for the proof.
\end{proof}
\end{prop}

\subsection{Elements of Sobolev calculus}\label{ss:SobolCal}
Now we want to concern certain notions related to Sobolev calculus on metric measure spaces.

As the extension of the classical energy appearing in \eqref{eq:OrigRes} to the singular context, in the standard way we adopt the so-called Cheeger energy, introduced by J. Cheeger in \cite{C99}. Among various essentially equivalent approaches to its definition, listed, for instance, in \cite{AILP24}, we choose the form proposed originally by M. Miranda in \cite{M03}, which is based on a specific relaxation procedure involving approximations by locally Lipschitz functions. We apply this construction only to real-valued functions, so, to distinct the notation from the case of maps valued in general metric spaces, we implement corresponding markers there.
\begin{defi}\label{def:ScChEnerg}
Let $p\in \mathbb{R}_{\geq 1}$, let $O\subseteq \X$ be an open set, let $u\in \Me(O,\R)$. Given an open set $O'\subseteq O$, \textbf{the scalar-valued Cheeger $p$-energy of $u$ on $O'$} is defined by
\begin{equation}
\mathrm{E}^{Sc}_p[u](O')\coloneqq \inf\limits_{(v_n)_{n\in \mathbb{N}}}\limi\limits_{n\rightarrow +\infty} \int\limits_{O'}\Big(\mathrm{lip}[v_n](x)\Big)^p\D\mathfrak{m}(x),
\end{equation}
where the infimum is taken over all sequences $(v_n)_{n\in \mathbb{N}}\subseteq\mathrm{LIP}^{Loc}(O')$ converging to $u$ in $\Leb^p_{Loc}(O')$ as $n\rightarrow +\infty$.
\end{defi}
\noindent With the data from the above definition in mind, we extend in the usual way the set function $\mathrm{E}^{Sc}_p[u](\bdot)$, which is so far defined only on the family of open subsets of $O $, to an arbitrary Borel set $V\subseteq O$ by letting
\begin{equation}\label{eq:EnCarExt}
\mathrm{E}^{Sc}_p[u](V)\coloneqq \inf\limits_{O'}\mathrm{E}^{Sc}_p[u](O'),
\end{equation}
where the infimum is taken over all open sets $O'\subseteq O$ with $V\subseteq O'$.

As it turns out, the set function obtained with the above procedure is a well-posed measure.
\begin{lem}\label{lem:ChEnerMeas}
Let $p\in \mathbb{R}_{\geq 1}$, let $O\subseteq \X$ be an open set, let $u\in\Me(O,\R)$. Then $\mathrm{E}^{Sc}_p[u](\bdot)$ is a Borel measure on $O$. Furthermore, given any open set $O'\subseteq O$, one has
\begin{equation}\label{eq:ChEnerMeas}
    \E_p^{Sc}[u](O')=\sup\limits_{V}\E_p^{Sc}[u](V),
\end{equation}
where the supremum is taken over all compact sets $V\subseteq O'$.
\begin{proof}
The first part is now a common fact, and its proof, for the case $p=1$, can be found in \cite[Theorem 3.4]{M03}. And the same arguments, with appropriate modifications, are also applicable to the case $p>1$.

To derive the second part, we can argue as follows. If the right-hand side of \eqref{eq:ChEnerMeas} is infinite, then there is nothing to prove. In the opposite case, we notice that $\mathrm{E}^{Sc}_p[u]\big(\{x\}\big)<+\infty$ for any $x\in O'$, which, together with \eqref{eq:EnCarExt}, implies that $\mathrm{E}_p[u](\bdot)$ is locally finite on $O'$. It then remains to use that locally finite, Borel measures on Polish spaces, which are known to include locally complete, separable metric spaces, satisfy a suitable form of the inner regularity property, which is exactly what we want. Thus, the equality in \eqref{eq:ChEnerMeas} is indeed always fulfilled.
\end{proof}
\end{lem}

We move to the definition of Cheeger energies for metric-valued maps. It relies on considering post-compositions of a given map with $1$-Lipschitz real-valued functions, which reduces everything to the scalar case. Originally, it was proposed by L. Ambrosio in \cite{A90} and by Y. G. Reshetnyak in \cite{R97}, and now it is recognized as the most appropriate one.
\begin{defi}\label{def:MetrCh}
Let $p\in \R_{\geq 1}$, let $(\Y, \dY)$ be a metric space, let $O\subseteq\X$ be an open set, let $f\in \Me(O,\Y)$. We define {\textbf{the Cheeger $p$-energy of $f$}}, which we denote by $\mathrm{E}_p[f]$, as the lowest element, in the sense of the lattice of non-negative Borel measures on $O$, among those $\nu$ of this type that meet the following:
\begin{equation}
\begin{gathered}
\mathrm{E}^{Sc}_p[\phi\circ f](V)\leq \nu(V)\\
\text{for any Borel set $V\subseteq O$ and any $\phi\in \mathrm{BLIP}_1(\Y)$}.
\end{gathered}
\end{equation}
\end{defi}
\noindent The above definition is well posed, since the considered lattice is known to be complete, and hence the stated object, which is basically just a lattice supremum due to \cref{lem:ChEnerMeas}, always exists. For a necessary information on the mentioned notions, we refer the reader to \cite{HKST15}.

\begin{rema}
    The construction appearing in \cref{def:ScChEnerg} can be easily seen applicable even to Banach-valued functions, not just real-valued ones. With this in mind, one may think about considering within \cref{def:MetrCh} post-compositions with bounded $1$-Lipschitz functions valued in arbitrary Banach spaces. This indeed gives rise to an object, a priori different, that can also be taken on the role of energy. And we claim that all the subsequent proofs can be carried out with this modified definition without any significant changes, as it was done in detail within \cite{O25}.

    Similarly, one may also want to deal with all $1$-Lipschitz functions, not just bounded ones, in \cref{def:MetrCh}. As opposed to \cite{O25}, this change is possible within our current exposition only partially. Loosely speaking, the issue here lies in that the post-composition of a metric-valued map with a general $1$-Lipschitz function may not be locally integrable, which is clearly not the case for a bounded $1$-Lipschitz function. And this circumstance really affects a part of the reasoning given in \cref{theo:LowBoundTh}, which should, in turn, be compared with \cite[Theorem 3.7]{O25}, where certain related things stand differently. Since all this is not actually relevant for us, we do not provide more details here.
\end{rema}

The following technical proposition, basically related to one straightforward property of lattice suprema, will be used later.
\begin{prop}\label{prop:LatSupExpr}
Let $p\in \R_{\geq 1}$, let $(\Y, \dY)$ be a metric space, let $O\subseteq\X$ be an open set, let $f\in \Me(O,\Y)$. Then one has
\begin{equation}
\label{eq:LatSupExpr}
\E_p[f](O)=\sup_{(O_j,\phi_j)_{j\in \mathbb{N}}} \sum\limits_{j\in \mathbb{N}} \mathrm{E}^{Sc}_p[\phi_j\circ f](O_j),
\end{equation}
where the supremum is taken over all families $(O_j,\phi_j)_{j\in \mathbb{N}}$, where, for each $j\in \mathbb{N}$, $O_j$ is an open subset of $O$ and $\phi_j$ is a function in $\mathrm{BLIP}_1(\Y)$, such that the family $(O_j)_{j\in \mathbb{N}}$ is pairwise disjoint.
\begin{proof}
The statement can be derived from some pretty standard measure-theoretic methods.
\end{proof}
\end{prop}

For the reader's convenience, we introduce the following objects, which, however, play a somewhat indirect role in the current exposition. Let $p\in \mathbb{R}_{\geq 1}$, let $(\Y, \dY)$ be a metric space, let $O\subseteq \X$ be an open set, let $f\in \Me(O,\Y)$. Given $x\in O$ and $r\in \mathbb{R}_{>0}$ with $\mathrm{B}(x,r)\subseteq O$, we put
\begin{gather}
\mathcal{A}_p[f](x,r)\coloneqq \begin{dcases}\frac{\mathrm{E}_p[f]\big(\mathrm{B}(x, r)\big)}{\mathfrak{m}\big(\mathrm{B}(x,r)\big)},\quad \mathfrak{m}\big(\mathrm{B}(x,r)\big)\in \mathbb{R}_{>0},\label{eq:AverRegOp}\\
0, \quad \text{otherwise},
\end{dcases}\\
\mathcal{R}_p[f](x,r)\coloneqq \tfrac{1}{3}\sum\limits_{k=0}^{+\infty}\big(\tfrac{2}{3}\big)^k\mathcal{A}_p[f]\big(x,\tfrac{r}{2^k}\big)\label{eq:RieszRegOp}.
\end{gather}
The first quantity above, being the averaging operator applied to $\E_p[f]$, is provided here only to define the second one, which can be viewed as a modified version of the so-called Riesz potential applied to $\mathrm{E}_p[f]$. The latter one will appear later within \cref{prop:PointIneq}, which, in turn, is presented without a detailed proof. So, for the specifics on how these very objects pop up in our context, see \cite[Lemma 3.3]{O25}.

For the sake of accuracy, we record the following technical proposition.
\begin{prop}
Let $p\in \R_{\geq 1}$, let $(\Y,\dY)$ be a metric space, let $O\subseteq \X$ be an open set, let $f\in \Me(O,\Y)$ be with $\E_p[f](O)<+\infty$. Let $V\Subset O$ be a Borel set, let $r\in \R_{>0}$ be such that $\B(V,r)\subseteq O$. Then the function $\mathcal{R}_p[f](\,\cdot\,, r)$ is Borel on $V$.
\begin{proof}
    For a possible proof, one may use the arguments provided in \cite[Section 3]{HKST15}.
\label{prop:RieszBorel}
\end{proof}
\end{prop}

Similar to \cite{O25}, it will be necessary for us to impose a certain extension property on open subsets of $\X$ suitably adapted to the metric-valued setting, for which we need the following auxiliary terminology. Given a metric space $(\Y, \dY)$, a set $S\subseteq \X$, a map $f\colon S\to\Y$, and a metric space $(\Z,\dZ)$, we say that a map $F\colon \X\to \Z$ is a $\Z$-extension of $f$ if there exists an isometric embedding $\iota\colon \Y\to \Z$ such that $F(x)=\iota\big(f(x)\big)$ for $\mathfrak{m}$-a.e. $x\in S$.

Using the notion introduced above, we provide the following definition.
\begin{defi}\label{def:StrExtDom}
Let $p\in \R_{\geq 1}$, let $(\Y,\dY)$ be a metric space, let $O\subseteq \X$ be an open set. We say that $O$ has {\textbf{the strong $p$-extension property with respect to $\Y$}} if any map $f\in \Leb^p(O, \Y)$ admits, for some metric space $(\Z,\dZ)$, a $\Z$-extension $F\in \Me(\X,\Z)$ such that
\begin{equation}
\label{eq:StrExtDom}
\lim\limits_{R\searrow 0} \mathrm{E}_p[F]\big(\mathrm{B}(O,R)\big)=\mathrm{E}_p[F](O).
\end{equation}
\end{defi}
\noindent For the sake of clarity, we note that it would change nothing above if we considered there only maps of finite energy.

\begin{rema}
The extension property presented is a metric-valued analog of the real-valued one used in \cite{LPZ22}. And it is a meaningful question of how they are related to each other depending on both the source and target spaces under consideration. We, however, leave this issue outside the scope and do not discuss it in more detail hereafter.
\end{rema}

\subsection{Doubling and Poincar{\'e} conditions}\label{ss:DoubPonCon}
In this section we discuss those structural assumptions on the underlying space that we use throughout the manuscript.

\begin{rema}
For the modern analysis on metric measure spaces, it is fairly common to require the ambient space to support some kinds of so-called doubling and Poincar{\'e} conditions, as such ones turn out to be simultaneously quite general and powerful. For more on their importance, we send the reader to \cite{HKST15}. And these conditions are exactly what we adopt for the current exposition. 

In our work, we wish to provide the most general possible versions of our statements, the fact of which affects the way we formulate the mentioned assumptions. In particular, we will deal with their appropriately ``localized'' variations, in a way similar to \cite{O25}. However, in contrast with \cite{O25}, we prefer here not to require the space to have full support, which forces us to modify correspondingly all the appearing definitions. Moreover, the conditions therein will be subordinate to a reference open set, as it enables us to cover a slightly broader context.
\end{rema}

We start by introducing the doubling property necessary for us. It will be used to prove the upper bound in \cref{theo:Res}.
\begin{defi}\label{def:DoubCon}
Let $\Omega\subseteq \X$ be an open set. We say that $\X$ is {\textbf{doubling around $\Omega$}} if there exists a constant $C\coloneqq C_{\mathrm{D}}[\Omega]$ and there is a radius $R\coloneqq R_{\mathrm{D}}[\Omega]\in \R_{>0}$ such that it holds, for any $x\in \B(\Omega,R)\cap \supp(\m)$ and any $r\in (0,R]$, that
\begin{equation}
\mathfrak{m}\big(\mathrm{B}(x,2r)\big)\leq C\mathfrak{m}\big(\mathrm{B}(x,r)\big).
\end{equation}
\end{defi}
\noindent We will refer to the inequality above as the doubling inequality.

Under the doubling condition presented above, the following properties of the Riesz potentials defined via \eqref{eq:RieszRegOp} turn out to be valid.
\begin{prop}\label{prop:DoubProp}
Let $\Omega\subseteq \X$ be an open set, suppose $\X$ is doubling around $\Omega$. Let $p\in \R_{\geq 1}$, let $(\Y, \dY)$ be a metric space, let $O\subseteq \X$ be an open set, let $f\in \Me(O,\Y)$. Let $E\Subset O$ be a $\m$-measurable set with $E\subseteq \Omega$, put $R\coloneqq \min\Big\{R_{\mathrm{D}}[\Omega], \dist\big(E, \X\backslash O\big)\Big\}$. Then the following assertions hold.
\begin{enumerate} 
\item[$\rm i)$] Given any $r\in (0,R]$, it holds, for $\m$-a.e. $x\in E$, that
\begin{equation}\label{eq:DoubProp1}
    \sup\limits_{r'\in \left[\frac{r}{2},r\right]}\mathcal{R}_p[f](x,r')\leq C_{\mathrm{D}}[\Omega]\mathcal{R}_p[f](x,r).
\end{equation}
\item[$\rm ii)$] If $\mathrm{E}_p[f](O)<+\infty$, then it holds, for any $r\in (0,R]$, that
\begin{equation}\label{eq:DoubProp2}
\int\limits_{E}\mathcal{R}_p[f](x,r)\D \mathfrak{m}(x)\leq C_{\mathrm{D}}[\Omega]\mathrm{E}_p[f]\big(\mathrm{B}(E,r)\big).
\end{equation}
\end{enumerate}
\begin{proof}
The inequality as in \eqref{eq:DoubProp1}, but with $\mathcal{R}_p[f]$ replaced by $\mathcal{A}_p[f]$, can be derived via the doubling inequality, which then leads easily to the desired estimate and hence proofs assertion $\rm i)$. In turn, assertion $\rm ii)$ can be obtained in a way similar to the one used in \cite[Lemma 3.11]{GT21}.
\end{proof}
\end{prop}

In addition to the pretty standard doubling condition from \cref{def:DoubCon}, we wish also to introduce one more related notion that pops up as a necessary requirement for the lower-bound part of \cref{theo:Res}.
\begin{defi}\label{def:StrDoub}
Let $\Omega\subseteq \X$ be an open set. We say that $\X$ is {\textbf{strongly doubling within $\Omega$}} if there exists a constant $C\coloneqq C_{\mathrm{SD}}[\Omega]\in \R_{\geq 1}$ and, for any compact set $V\subseteq \Omega$, there is a radius $R\coloneqq R_{\mathrm{SD}}[\Omega,V]$ such that it holds, for any $x\in \B(V,R)\cap \supp(\m)$ and any $r\in (0,R]$, that
\begin{equation}
    \m\big(\B(x,2r)\big)\leq C\m\Big(\A\big(x,r,\tfrac{3r}{4}\big)\Big).
\end{equation}
\end{defi}
\noindent The inequality above will be referred to as the strong doubling one.

\begin{rema}
    In contrast to the usual doubling inequality, the strong doubling one allows us to control, for a given ball, the measure of the doubled ball by the measure not only of the ball itself, but also of a certain spherical shell in it. Whereas this property is a priori stronger than the former one, both turn to be equivalent under certain mild conditions on the ambient space, as indicated in \cref{prop:DoubStrDoub} below. We should also note that the factor $\tfrac{3}{4}$ is essentially arbitrary, yet its change will affect some of the universal constants that appear later on.
\end{rema}

The following properties implied by the strong doubling condition will be used later.
\begin{prop}\label{prop:StrDoubProp}
    Let $\Omega \subseteq \X$ be an open set, suppose $\X$ is strongly doubling within $\Omega$. Let $O\subseteq \Omega$ be an open set. Then the following assertions hold.
    \begin{enumerate}
        \item[$\rm i)$] Let $p\in \R_{\geq 1}$, let $u\in \Leb^p_{Loc}(O)$. Then $\m$-a.e. point $x\in O$ is a $p$-Lebesgue point of $u$.
        \item[$\rm ii)$] Let $u\in \Leb^1(O)$, let $O'\Subset O$ be an open set. Put $R\coloneqq \dist\big(O',\X\backslash O\big)$. Then the functions $\langle u\rangle_{\B(\,\cdot\,,r)}$, $r\in (0,R)$, converge to $u$ in $\Leb^1_{Loc}(O')$ as $r\searrow 0$.
    \end{enumerate}
    \begin{proof}
        Assertion $\rm i)$ is the standard improved version of the classical Lebesgue differentiation theorem, which is known to hold under doubling conditions even much weaker than ours. For the proof, see \cite[Section 4.1]{HKST15} for instance. In turn, the proof of assertion $\rm ii)$ can be carried out easily via \cite[Lemma 3.11]{GT21}.
    \end{proof}
\end{prop}
\noindent Clearly, for the just given statement the very strong doubling inequality is excessive, and the standard one completely suffices here as well.

The following proposition provides a link between the two introduced doubling conditions.
\begin{prop}\label{prop:DoubStrDoub}
    Let $\Omega\subseteq \X$ be a connected, open set with $\Omega\subseteq \supp(\m)$, suppose $\X$ is doubling around $\Omega$. Then $\X$ is strongly doubling within $\Omega$ and, moreover, the constant $C_{\mathrm{SD}}[\Omega]$ can be chosen so that it depends only on $C_{\mathrm{D}}[\Omega]$. Furthermore, one has $\m\big(\{x\}\big)=0$ for every $x\in \Omega$.
    \begin{proof}
        The first statement easily follows from the so-called reverse doubling inequality, the validity of which in our situation can be verified via \cite[Remark 8.1.15]{HKST15}. Exactly the same inequality can be used to verify the second statement.
    \end{proof}
\end{prop}
\noindent The second part of the above proposition shows that, under the corresponding assumptions, the diagonal of our space is of measure
zero, the fact of which should be kept in mind when deriving \cref{theo:Res} from the subsequent statements.

Now we introduce a Poincar{\'e} condition that will be considered in what follows within the proof of the upper-bound part of \cref{theo:Res}.
\begin{defi}\label{def:PoinIneq}
Let $\Omega\subseteq \X$ be an open set, let $p\in \mathbb{R}_{\geq 1}$. We say that $\X$ supports {\rm\textbf{a $p$-Poincar{\'e} inequality around $\Omega$}} if there exists a constant $C\coloneqq C_{\mathrm{P}}[\Omega]\in \mathbb{R}_{>0}$ and there are a parameter $\lambda\coloneqq\lambda_{\mathrm{P}}[\Omega]\in \mathbb{R}_{\geq 1}$ and a radius $R\coloneqq R_{\mathrm{P}}[\Omega]\in \mathbb{R}_{>0}$ such that, for any $x\in \B(\Omega,R)\cap \supp(\m)$ and any $r\in (0,R]$ with $\m\big(\B(x,\lambda r)\big)<+\infty$, one has
\begin{equation}
\fint\limits_{\mathrm{B}(x,r)}\Big|u(x')-\langle u\rangle_{\mathrm{B}(x,r)}\Big|\D\mathfrak{m}(x')\leq C r \Biggg(\fint\limits_{\mathrm{B}(x,\lambda r)}\Big(\mathrm{lip}[u](x')\Big)^p\D \mathfrak{m}(x')\Biggg)^\frac{1}{p}
\end{equation}
for every Lipschitz function $u\colon \X\to \mathbb{R}$.
\end{defi}

As a convenient unifying notion, the very one used in \cref{theo:Res}, we introduce the following.
\begin{defi}\label{def:PoinSpace}
    Let $p\in \R_{\geq 1}$. We say that $\X$ is \textbf{a $p$-Poincar{\'e} space} if it is connected, has full support, is doubling around $\X$, and supports a $p$-Poincar{\'e} inequality around $\X$. 
\end{defi}

\begin{rema}
    As it follows immediately from \cref{prop:DoubStrDoub}, if $\X$ is a Poincar{\'e} space, then it is also strongly doubling within $\X$. Thus, under the former condition, all the properties listed in the current section become perfectly valid.
\end{rema}

\subsection{Requirements for mollifiers}\label{ss:ReqMol}
Here we indicate the assumptions on mollifiers that are exploited subsequently.

Following the style from \cite{O25}, we prefer to formulate the conditions of mollifiers in a quite cumbersome way, exactly as they emerge within our reasoning, for which we shall also introduce some special notation. And later, namely in \cref{ss:SimplAssMol}, we will provide their highly simplified versions, which have already a much closer form to what one can see in \eqref{eq:OrigNormalCond} and \eqref{eq:OrigLocalCond}.

Within the context as below, the notation like $\tau_{\bullet}$ is reserved for sequences $(\tau_k)_{k\in \mathbb{N}_0}$ of intervals in $\R$, while the notation like $\rho_{\bullet}$ is reserved for families $(\rho_{\delta})_{\delta\in (0,1)}$ of functions $\X\times \X\to \overline{\R}_{\geq 0}$.

Let $r\in \mathbb{R}_{>0}$. By $\P_r^L$ we denote the family of all sequences $\tau_{\bullet}$ with the following properties: 
\begin{enumerate}
\item[$\bullet$] for each $k\in \mathbb{N}_0$, $\tau_k$ is a non-empty interval in $(0,r]$ with $2\inf(\tau_k)\leq \sup(\tau_k)$;
    \item[$\bullet$] for all distinct $k,k'\in \mathbb{N}_0$, it holds that
    \begin{equation}
        \tau_k\cap \tau_{k'}=\varnothing.
    \end{equation}
\end{enumerate}
Given $\rho\in \Me_+\big(\X\times \X\big)$ and a $\m$-measurable set $E\subseteq \X$, we put first
\begin{equation}\label{eq:LowIntMolPart}
    \mathcal{I}^L_{E,r}[\rho,\tau_{\bullet}]\coloneqq \sum\limits_{k\in \mathbb{N}_0}\einf\limits_{x\in E} \bigg(\m\big(\A(x,\tau_k)\big)\einf\limits_{x'\in \A(x,\tau_k)}\rho(x,x')\bigg)
\end{equation}
for any $\tau_{\bullet}\in \P^L_r$ and then \begin{equation}\label{eq:LowIntMolFul}
    \mathcal{I}^L_{E,r}[\rho]\coloneqq \sup\limits_{\tau_{\bullet}\in \P^L_r} \mathcal{I}^L_{E,r}\big[\rho,\tau_{\bullet}\big].
    \end{equation}
    
Now, given $\rho_{\bullet}\subseteq \Me_+\big(\X\times \X\big)$ and a compact set $V\subseteq \X$, we define the quantity
\begin{equation}\label{eq:LowIntMolFamily}
        \mathcal{I}^L_V[\rho_{\bullet}]\coloneqq \sup\limits_{E}\lims\limits_{r\searrow 0}\limi\limits_{\delta\searrow 0}\mathcal{I}^L_{E,r}[\rho_{\delta}],
    \end{equation}
    where the supremum is taken over all $\m$-measurable sets $E\subseteq \X$ with $V\Subset E$.

The condition necessary for the proof of the lower bound in \cref{theo:Res} is the following.
\begin{defi}\label{def:LowAdmFam}
    Let $\Omega\subseteq \X$ be an open set, let $\rho_{\bullet}\subseteq \Me_+\big(\X\times \X\big)$. We say that $\rho_{\bullet}$ is {\textbf{lower-admissible for $\Omega$}} if there exists a constant $C\coloneqq C_{\mathrm{M}}^L[\Omega,\rho_{\bullet}]\in \R_{>0}$ such that
    \begin{equation}\label{eq:LowAdmFam}
    \inf\limits_{V}\mathcal{I}^L_V[\rho_{\bullet}]\geq C,
\end{equation}
where the infimum is taken over all compact sets $V\subseteq \Omega$.
\end{defi}

Let $r\in \R_{>0}$. By $\P^U_r$ we denote the family of all sequences $\tau_{\bullet}$ with the following properties:
\begin{enumerate}\label{eq:LowAdmIneq}
    \item[$\bullet$] for each $k\in \mathbb{N}_0$, $\tau_k$ is a non-empty interval in $(0,r]$ with $\sup(\tau_k)\leq 2\inf(\tau_k)$;
    \item[$\bullet$] it holds that
    \begin{equation}
        (0,r]\subseteq \bigcup\limits_{k\in \mathbb{N}_0}\tau_k.
    \end{equation}
\end{enumerate}
Given $\rho\in \Me_+\big(\X\times \X\big)$ and $S\subseteq \X$, we put first
\begin{equation}\label{eq:UpIntMolPart}
    \mathcal{I}^U_{S,r}[\rho, \tau_{\bullet}]\coloneqq  \sum\limits_{k\in \mathbb{N}_0} \esup\limits_{x\in \B\bg(S,\sup(\tau_k)\bg)} \int\limits_{\A(x,\tau_k)}\Big(\rho(x,x')+\rho(x',x)\Big)\D \m(x')
\end{equation}
for any $\tau_{\bullet}\in \P_r^U$ and then
\begin{equation}\label{eq:UpIntMolFul}
    \mathcal{I}^U_{S,r}[\rho]\coloneqq \inf\limits_{\tau_{\bullet}\in \P_r^U} \mathcal{I}^U_{S,r}[\rho, \tau_{\bullet}].
\end{equation}

Now, given $\rho_{\bullet}\subseteq \Me_+\big(\X\times \X\big)$ and $S\subseteq \X$, we put
\begin{equation}\label{eq:UpIntMolFamily}
    \mathcal{I}^U_S[\rho_{\bullet}]\coloneqq \limi\limits_{r\searrow 0}\lims\limits_{\delta \searrow 0} \mathcal{I}^U_{S,r}[\rho_{\delta}].
\end{equation}

For the proof of the upper bound in \cref{theo:Res}, we will use the following condition.
\begin{defi}\label{def:UpAdmFam}
    Let $\Omega\subseteq \X$ be an open set, let $\rho_{\bullet}\subseteq \Me_+\big(\X\times \X\big)$. We say that $\rho_{\bullet}$ is {\textbf{upper-admissible for $\Omega$}} if there exists a constant $C\coloneqq C^U_{\mathrm{M}}[\Omega,\rho_{\bullet}]\in \R_{>0}$ such that
    \begin{equation}\label{eq:UpAdmFam}
        \mathcal{I}^U_{\Omega}[\rho_{\bullet}]\leq C.
    \end{equation}
\end{defi}

\begin{rema}
    Let us provide some brief clarification for the reader on what is introduced above. The collections $\mathcal{P}^L_r$ and $\mathcal{P}^U_r$, with $r\in \R_{>0}$, consist of special ``partitions'' of $(0,r]$ adapted to the subsequent proofs of the lower and upper bounds in \cref{theo:Res}, respectively. More explicitly, partitions in each of the collections are required to be, respectively, not too dense and not too sparse in a quantitative way, and this will play its role later. In turn, the objects appearing in \eqref{eq:LowIntMolPart} and \eqref{eq:UpIntMolPart} can be viewed as ``integral sums'' of a given mollifier function with respect to corresponding partitions. These quantities should also be compared with the integral from \eqref{eq:OrigNormalCond}, as the former ones express basically the same thing, but in a more complicated way. So one can see direct parallels of the conditions in \cref{def:LowAdmFam} and \cref{def:UpAdmFam} with the original one.
\end{rema}

\begin{rema}
    It is worth noticing that the conditions in \cref{def:LowAdmFam} and \cref{def:UpAdmFam} do not involve any dependence on $p$, as it is the case for the original BBM assumptions given by \eqref{eq:OrigNormalCond}. And whereas this feature is in line with what took place within \cite{LPZ22} and \cite{O25} for the upper-bound condition, it is somewhat novel for the lower-bound condition.
\end{rema}

In addition to the previous conditions, which impose restrictions on the behavior of mollifiers near the diagonal, we also need the following one, which constrains the behavior of such at infinity.
\begin{defi}\label{def:MolDecCond}
Let $\Omega\subseteq \X$ be an open set, let $p\in \R_{\geq 1}$, let $\rho_{\bullet}\subseteq \Me_+\big(\X\times \X\big)$. We say that $\rho_{\bullet}$ satisfies {\textbf{the $p$-decay condition for $\Omega$}} if
\begin{equation}\label{eq:MolDecCond}
\lim\limits_{r\searrow 0}\lims\limits_{\delta \searrow 0} \esup\limits_{x\in \Omega} \int\limits_{\Omega\cap \mathrm{A}(x,\infty,r)} \frac{\rho_{\delta}(x,x')+\rho_{\delta}(x',x)}{\big(\dX(x,x')\big)^p} \D \mathfrak{m}(x')=0.
\end{equation}
\end{defi}
\noindent As with \cref{def:LowAdmFam} and \cref{def:UpAdmFam}, there should be noticed a resemblance of the condition above with what is given by \eqref{eq:OrigLocalCond}.

At last, we can formulate a convenient notion unifying all the previous ones, which exactly appears in \cref{theo:Res}.
\begin{defi}\label{def:AdmFam}
    Let $p\in \R_{\geq 1}$, let $\rho_{\bullet}\subseteq \Me_+\big(\X\times \X\big)$. We say that $\rho_{\bullet}$ is {\textbf{$p$-admissible}} if $\rho_{\bullet}$ is upper-admissible and lower-admissible and satisfies the $p$-decay condition, all for $\X$.
\end{defi}

\section{Proofs}\label{ss:Proofs}

\subsection{Lower bound}\label{ss:LowBound}
In this section we prove the lower-bound estimate from \cref{theo:Res}.

A necessary technical tool for us is certain Lipschitz partitions of unity, quite similar to those that appear in \cite[Lemma 2.4]{GT21}. However, as our formulation differs to some extent from the one therein, we prefer to provide an independent proof for the reader's convenience.
\begin{prop}\label{prop:PartUnit}
Let $\Omega\subseteq \X$ be an open set, suppose $\X$ is strongly doubling within $\Omega$. Then there exists a constant $C=C(\Omega)\in \mathbb{R}_{\geq 1}$ depending only on $C_{\mathrm{SD}}[\Omega]$ such that the following holds. Let $V\subseteq \Omega$ be a compact set. Then there is a radius $R=R(\Omega,V)\in \mathbb{R}_{>0}$ such that, putting $S\coloneqq \B(V,R)\cap \supp(\m)$, one can find, for any $r\in (0,R]$, a data set, parametrized by some at most countable index set $I$, consisting of points $x_i\in S$, $i\in I$, and of functions $\varphi_i\colon \X\to  [0,1]$, $i\in I$, with the following properties:
\begin{itemize}
\item[$\rm i)$] for every $x\in \X$ one has $\sum\limits_{i\in I} \chi_{\mathrm{B}(x_i,2r)}(x)\leq C$;
\item[$\rm ii)$] for each $i\in I$ it holds that $\Lip[\varphi_i]\leq \tfrac{C}{r}$ and that $\supp(\varphi_i)\subseteq \B(x_i,r)$;
\item[$\rm iii)$] for every $x\in \B\big(S,\tfrac{r}{16}\big)$ one has $\sum\limits_{i\in I} \varphi_i(x)=1$.
\end{itemize}
    \begin{proof}
Put for brevity $C_{\mathrm{SD}}\coloneqq C_{\mathrm{SD}}[\Omega]$. Define $R\coloneqq \tfrac{1}{16} R_{\mathrm{SD}}[\Omega,V]$ and fix $r\in (0,R]$.
    
        We clearly can find a maximal $\tfrac{r}{4}$-separated set $(x_i)_{i\in I}$ in $S$, where $I$ is some at most countable index set. That is, it holds that $\dX(x_i,x_{i'})>\tfrac{r}{4}$ for all distinct $i,i'\in I$ and that for any $x\in S$ there is $i\in I$ with $x\in \B\big(x_i,\tfrac{r}{4}\big)$. Exactly the points $x_i$, $i\in I$, constitute the declared data.
        
        We want now to estimate how much the balls $\B(x_i,2r)$, $i\in I$, overlap. Fix any $x\in \X$ and put
        \begin{equation}
            I_x\coloneqq \Big\{i\in I\bigm| \B(x,2r)\cap \B(x_i,2r)\neq \varnothing\Big\}.
        \end{equation}
        We claim that
        \begin{equation}\label{eq:PartUnit1}
            \card(I_x)\leq( C_{\mathrm{SD}})^7.
        \end{equation}
        Indeed, if $I_x$ is empty, then there is nothing to prove. Thus we can assume that there is $i_0\in I_x$. This implies that $\dX(x_{i_0},x)\leq 4r$, whence we have $\dX(x_{i_0},x_i)\leq 8r$ for any $i\in I_x$. This ensures us that the balls $\B\big(x_i,\tfrac{r}{8}\big)$, $i\in I_x$, are contained in $\B\big(x_{i_0},16r\big)$. Since they are also disjoint, we can write, with the use of the doubling inequality, that
        \begin{equation}
        \begin{gathered}
            \m\Big(\B\big(x_{i_0},16r\big)\Big)\geq \sum\limits_{i\in I_x} \m\Big(\B\big(x_i,\tfrac{r}{8}\big)\Big)\geq \\
            \geq \frac{1}{(C_{\mathrm{SD}})^7}\sum\limits_{i\in I_x} \m\Big(\B\big(x_i,32r\big)\Big)\geq \frac{\card(I_x)}{(C_{\mathrm{SD}})^7}\m\Big(\B\big(x_{i_0},16r\big)\Big).
            \end{gathered}
        \end{equation}
        Since $x_{i_0}\in \supp(\m)$, we get the estimate in \eqref{eq:PartUnit1}.

        Now we move to defining corresponding functions for the data. To begin with, we pick $i\in I$ and define a function $\psi_i\colon \X\to \R_{\geq 0}$ by
        \begin{equation}\label{eq:PartUnit2}
            \psi_i(x)\coloneqq \max\Big\{\tfrac{r}{2}-\dX(x_i,x),0\Big\}.
        \end{equation}
        It is obvious that $\Lip[\psi_i]\leq 1$, that $\supp(\psi_i)\subseteq \B(x_i,r)$, and that $\sup(\psi_i)\leq \tfrac{r}{2}$.

        Put then $S'\coloneqq \B\big(S,\tfrac{r}{16}\big)$ and make the following simple observation: if $x\in S'$, then there is $i\in I$ with $\dX(x,x_i)\leq \tfrac{3r}{8}$. With this, from \eqref{eq:PartUnit2} it follows directly that
        \begin{equation}\label{eq:PartUnit3}
            \inf\limits_{x\in S'}\sum\limits_{i\in I} \psi_i(x)\geq \frac{r}{8}.
        \end{equation}
        Define then a function $\Psi\colon \X\to \R_{\geq 0}$ by
        \begin{gather}\label{eq:PartUnit4}
            \Psi\coloneqq \max\Bigg\{\sum\limits_{i\in I} \psi_{i},\frac{r}{8}\Bigg\},
       \end{gather}
        which is well posed, since the sum appearing above is locally finite due to what is established previously. So, we see that
        \begin{equation}
            \inf(\Psi)\geq \frac{r}{8}.
        \end{equation}
        Moreover, given any $i\in I$, we can easily deduce, again from the previous estimates, that
        \begin{equation}
            \Lip\Big[\Psi\big|_{\B(x_i,2r)}\Big]\leq \sum\limits_{i'\in I} \Lip\Big[\psi_{i'}\big|_{\B(x_i,2r)}\Big]\leq (C_{\mathrm{SD}})^7,
        \end{equation}
        whence it also follows, with the use of assertion $\rm ii)$ from \cref{prop:GlLipConstProp}, that
        \begin{equation}
            \Lip\bigg[\frac{1}{\Psi}\Big|_{\B(x_i,2r)}\bigg]\leq \frac{64(C_{\mathrm{SD}})^7}{r^2}.
        \end{equation}

        Now, for each $i\in I$, we define a function $\varphi_i\colon \X\to [0,1]$ by
        \begin{equation}
            \varphi_i\coloneqq \frac{\psi_i}{\Psi},
        \end{equation}
        which clearly grants that $\supp(\varphi_i)\subseteq \B(x_i,r)$. Also, it follows from \eqref{eq:PartUnit3} and \eqref{eq:PartUnit4} that property $\rm iii)$ indeed takes place. We are left to estimate the constants $\Lip[\varphi_i]$, $i\in I$.

        Fix $i\in I$ and pick arbitrary $x,x'\in \X$. If $x,x'\not\in \supp(\varphi_i)$, then $\Q_{\varphi_i}(x,x')=0$. In turn, if $x,x'\in \B(x_i,2r)$, then we can write, with the use of assertions $\rm i)$ and $\rm ii)$ from \cref{prop:GlLipConstProp}, that
        \begin{equation}
        \begin{gathered}
            \Q_{\varphi_i}(x,x')\leq \Lip\Big[\varphi_i\big|_{\B(x_i,2r)}\Big] \leq \\
            \leq \frac{\Lip[\psi_i]}{\inf(\Psi)}+\sup(\psi_i) \Lip\bigg[\frac{1}{\Psi}\Big|_{\B(x_i,2r)}\bigg]\leq \frac{128 (C_{\mathrm{SD}})^7}{r}.
            \end{gathered}
        \end{equation}
        For the last case, we can assume without loss of generality that $x\in \supp(\varphi_i)$ and $x'\not\in \B(x_i,2r)$. Then $\dX(x,x')\geq r$, which gives us the estimate
        \begin{equation}
        \begin{gathered}
            \Q_{\varphi_i}(x,x')\leq \frac{\sup(\psi_i)}{r\inf(\Psi)}\leq \frac{4}{r}.
            \end{gathered}
        \end{equation}

        It remains to put $C\coloneqq 128 (C_{\mathrm{SD}})^7$ and notice that properties $\rm i)$ and $\rm ii)$ are thus valid. The proof is fully complete.
    \end{proof}
\end{prop}
\noindent We note that the validity of exactly the strong doubling inequality is clearly redundant for the statement above to hold. It suffices to have the standard doubling inequality only.

In the lemma below, we demonstrate the approximation procedure that will be used to obtain the necessary bound for the energy. While the corresponding technique is still based on discrete convolutions, as well as in \cite[Lemma 3.4]{O25}, there is one pivotal difference: instead of averaging a given function over balls, we deal with averages over annuli. As one can trace from the subsequent reasoning, this very modification allows us to extend the class of mollifiers that can be applied to the final result. 
\begin{lem}\label{lem:ApproxSeq}
    Let $\Omega\subseteq \X$ be an open set, suppose $\X$ is strongly doubling within $\Omega$. Let $p\in \R_{\geq 1}$. Then there exists a constant $C=C(\Omega)\in \R_{\geq 1}$ depending only on $p$, $C_{\mathrm{SD}}[\Omega]$ such that the following holds. Let $V\subseteq \Omega$ be a compact set. Then there is a radius $R=R(\Omega,V)\in \R_{>0}$ with the property that, given any $r\in (0,R]$, any $S\subseteq \B(V,R)$, and any $u\in \Leb^1\big(\B(S,r)\big)$, one can construct $u^r\in \LIP^{Loc}(\X)$ such that, for $\m$-a.e. $x\in S$, one has
    \begin{gather}
\Big|u^r(x)-u(x)\Big|^p\leq C \fint\limits_{\B(x,r)}\big(\dX_u(x,x')\big)^p\D\m(x'),\label{eq:ApproxSeq0}\\
 \Big(\mathrm{lip}\big[u^r\big](x)\Big)^p\leq \frac{C}{\m\big(\B(x,r)\big)}\int\limits_{\A\big(x,\left(\frac{r}{2},r\right)\big)}\big(\Q_u(x,x')\big)^p\D\m(x').\label{eq:ApproxSeq00}
\end{gather}
    \begin{proof}
        Let $C'=C'(\Omega)$ denote the constant and let $R'=R'(\Omega,V)$ denote the radius guaranteed by \cref{prop:PartUnit}. Find then $R_0\in \R_{>0}$ with $\m\big(\B(V,R_0)\big)<+\infty$. Now we define $R\coloneqq  \min\big\{R', R_{\mathrm{SD}}[\Omega,V], R_0\big\}$, after which we put $S_0\coloneqq \B(V,R)\cap \supp(\m)$ and fix any $r\in (0,R]$. Putting $r'\coloneqq \tfrac{r}{32}$, we apply \cref{prop:PartUnit} with respect to $S_0$ and $r'$, which gives us families $(x_i)_{i\in I}$ and $(\varphi_i)_{i\in I}$ with the corresponding properties. For each $i\in I$, we then also put $B_i\coloneqq \B\big(x_i,2r'\big)$ and $A_i\coloneqq \A\Big(x_i,\big(18r',30r'\big)\Big)$.

Now we fix arbitrary $S\subseteq \B(V,R)$ and $u\in \Leb^1\big(\B(S,r)\big)$. Define then a function $u^r\colon \X\to \R$ by
\begin{equation}\label{eq:ApproxSeq1}
u^r(x)\coloneqq\sum\limits_{i\in I} \varphi_i(x)\langle u\rangle_{A_i}.
\end{equation}
By the properties from \cref{prop:PartUnit}, it follows that the right-sided expression above is a locally finite sum of Lipschitz functions and hence defines a well-posed locally Lipschitz function.

Fix arbitrary $x\in S\cap \supp(\m)$. Note that $x\in S_0$. Also, by the strong doubling inequality, we easily have $\m(A_i)\in \R_{>0}$ for each $i\in I$.

The definition in \eqref{eq:ApproxSeq1} together with the properties from \cref{prop:PartUnit} imply that
\begin{equation}\label{eq:ApproxSeq2}
\begin{gathered}
\Big|u^r(x)-u(x)\Big|=\Bigg|\sum\limits_{i\in I}\varphi_i(x) \Big(\langle u\rangle_{A_i}-u(x)\Big)\Bigg| \leq \\
 \leq \sum\limits_{i\in I}\varphi_i(x) \Big|\langle u\rangle_{A_i}-u(x)\Big| \leq C' \sup\limits_{i\in I}\bigg(\chi_{B_i}(x)\Big|\langle u\rangle_{A_i}-u(x)\Big|\bigg).
\end{gathered}
\end{equation}
Pick then any $x'\in \B\big(x,\tfrac{r'}{16}\big)$ with $x'\neq x$. Carrying out similar computations as above and using again the properties from \cref{prop:PartUnit}, we can get
\begin{equation}\label{eq:ApproxSeq3}
\begin{gathered}
\Big|u^r(x')-u^r(x)\Big|
=\Bigg|\sum\limits_{i\in I}\Big(\varphi_i(x')-\varphi_i(x)\Big)\Big(\langle u\rangle_{A_i}-u(x)\Big)\Bigg|\leq \\
\leq  \frac{(C')^2}{r'}\mathsf{d}(x,x')\sup\limits_{i\in I} \bigg(\chi_{B_i}(x)\Big|\langle u\rangle_{A_i}-u(x)\Big|\bigg).
\end{gathered}
\end{equation}
Now we notice that for any $i\in I$ with $\dX(x,x_i)\leq 2r'$ it holds that
\begin{equation}\label{eq:ApproxSeq4}
\begin{gathered}
    \Big|\langle u\rangle_{A_i}-u(x)\Big|^p\leq \fint\limits_{A_i} \big(\Q_u(x,x'')\big)^p \D\m(x'')\leq\\
    \leq \frac{1}{\m\big(A_x''\big)}\int\limits_{A_x'}\big(\Q_u(x,x'')\big)^p \D\m(x''),
    \end{gathered}
\end{equation}
where we used Jensen's inequality and put $A_x'=\A\Big(x,\big(16r',32r'\big)\Big)$ and $A_x''\coloneqq \A\Big(x,\big(20r',28r'\big)\Big)$. The last step is to apply the strong doubling inequality and write
\begin{equation}\label{eq:ApproxSeq5}
    \m\Big(\B\big(x,32r'\big)\Big)\leq \m\Big(\B\big(x,54r'\big)\Big)\leq C_{\mathrm{SD}}[\Omega]\m\big(A_x''\big).
\end{equation}

Finally, we define $C \coloneqq C_{\mathrm{SD}}[\Omega]\big(32(C')^2\big)^p$ and recall that $r=32r'$. Once we combine \eqref{eq:ApproxSeq2}, \eqref{eq:ApproxSeq4}, and \eqref{eq:ApproxSeq5}, we can come easily to the estimate in \eqref{eq:ApproxSeq0}. In turn, if we pass to the upper limit as $x'\to x$ in \eqref{eq:ApproxSeq3} and combine it with \eqref{eq:ApproxSeq4} and \eqref{eq:ApproxSeq5}, we get the estimate in \eqref{eq:ApproxSeq00}. As $\m\big(S\backslash \supp(\m)\big)=0$, the proof is thus complete.
    \end{proof}
\end{lem}

\begin{rema}
    We note that in the previous lemma, all the same can be perfectly done for Banach-valued functions as well without any relevant modifications. In light of this, there can be achieved stronger variations of the remaining statements in the section, with a slightly different notion of energy involved for metric-valued maps, which is described in the remark under \cref{def:MetrCh}. For the specifics on this, we send the reader to \cite{O25}, where such a procedure was performed.
\end{rema}

The following intermediate statement shows how one can apply the approximants from the previous lemma in order to obtain a suitable estimate from above for the energy of a given function. For the proof, we recall \cref{prop:GenLebTh}, \cref{def:ScChEnerg}, \cref{prop:StrDoubProp}, as well as the notation preceding \cref{def:LowAdmFam}.
\begin{lem}\label{lem:LowBoundLem}
    Let $\Omega\subseteq \X$ be an open set, suppose $\X$ is strongly doubling within $\Omega$. Let $p\in \R_{\geq 1}$, let $\rho_{\bullet}\subseteq \Me_+\big(\X\times \X\big)$. Let $O\subseteq \X$ be an open set, let $u\in \Me(O,\R)$, let $O'\Subset O$ be an open set with $u\in \Leb^p_{Loc}(O')$. Then, for any compact set $V\subseteq O'\cap \Omega$, one has
    \begin{equation}\label{eq:LowBoundLem0}
    \begin{gathered}
        \frac{1}{C}\mathcal{I}^L_V[\rho_{\bullet}]\E^{Sc}_p[u](V)\leq\\
        \leq \lim\limits_{r\searrow 0} \limi\limits_{\delta\searrow 0}\int\limits_{O'}\int\limits_{\B(x,r)} \big(\Q_u(x,x')\big)^p \rho_{\delta}(x,x')\D \m(x')\D \m(x),
        \end{gathered}
    \end{equation}
    where $C=C(\Omega)$ is the constant guaranteed by \cref{lem:ApproxSeq}.
    \begin{proof}
        Fix a compact set $V\subseteq O'\cap \Omega$. Due to all the assumptions, we can find $R_0\in \R_{>0}$ with $\B(V,R_0)\Subset O'\cap \Omega$ such that $\m\big(\B(V,R_0)\big)<+\infty$ and $u\in \Leb^p\big(\B(V,R_0)\big)$. Let also $R'=R'(\Omega,V)$ denote the radius guaranteed by \cref{lem:ApproxSeq}. Then we put $R\coloneqq \tfrac{1}{2}\min\big\{R_0,R', R_{\mathrm{SD}}[\Omega,V]\big\}$.

        Fix a $\m$-measurable set $E\subseteq \X$ with $V\Subset E$ and pick then an open set $O''\subseteq \B(V,R)\cap E$ with $V\subseteq O''$. For every $r\in (0,R]$, let $u^r$ denote the function in $\LIP^{Loc}(\X)$ from \cref{lem:ApproxSeq} constructed with respect to $O''$, $r$, and $u$.

        Fix $r\in (0,R]$ and pick then any $\tau_{\bullet}\in \P^L_r$. For each $k\in \mathbb{N}_0$ put $t_k\coloneqq \sup(\tau_k)$. From the definition of $\P^L_r$, we know that $\big(\tfrac{t_k}{2},t_k\big)\subseteq \tau_k\subseteq (0,t_k]$ for each $k\in \mathbb{N}_0$, as well as that the intervals $\tau_k$, $k\in \mathbb{N}_0$, are pairwise disjoint. Fix also $\delta\in (0,1)$.
        
        With the use of \eqref{eq:ApproxSeq00} from \cref{lem:ApproxSeq}, we can write for $\m$-a.e. $x\in O''$ that
        \begin{equation}
            \begin{gathered}
                C\int\limits_{\B(x,r)} \big(\Q_u(x,x')\big)^p \rho_{\delta}(x,x')\D \m(x')\D \m(x)\geq \\
                \geq \sum\limits_{k\in \mathbb{N}_0} \m\big(\B(x,t_k)\big)\bigg(\einf\limits_{x'\in \A(x,\tau_k)} \rho_{\delta}(x,x')\bigg)\Big(\lip\big[u^{t_k}\big](x)\Big)^p\geq \\
                \geq \sum\limits_{k\in \mathbb{N}_0} \einf\limits_{x''\in E}\Bigg(\m\big(\A(x'',\tau_k)\big)\einf\limits_{x'\in \A(x'',\tau_k)} \rho_{\delta}(x'',x')\Bigg) \Big(\lip\big[u^{t_k}\big](x)\Big)^p.
            \end{gathered}
        \end{equation}
        Integrating the above estimate over $O''$, carrying out simple transformations, and using the notation given in \eqref{eq:LowIntMolPart}, we get
        \begin{equation}
            \begin{gathered}
                C\int\limits_{O'}\int\limits_{\B(x,r)} \big(\Q_u(x,x')\big)^p \rho_{\delta}(x,x')\D \m(x')\D \m(x)\geq \\
                \geq \mathcal{I}^L_{E,r}[\rho_{\delta},\tau_{\bullet}]\inf\limits_{k\in \mathbb{N}_0} \int\limits_{O''}\Big(\lip\big[u^{t_k}\big](x)\Big)^p\D \m(x).
            \end{gathered}
        \end{equation}
        Now we maximize the above inequality over $\tau_{\bullet}\in \P^L_r$, pass then there to the lower limit as $\delta\searrow 0$ and to the upper limit as $r\searrow 0$. All this, together with \eqref{eq:LowIntMolFul}, leads to the existence of $(r_n)_{n\in \mathbb{N}}\subseteq (0,R]$ with $\lims\limits_{n\to +\infty} r_n =0$ such that
         \begin{equation}
            \begin{gathered}
                C\lim\limits_{r\searrow 0}\limi\limits_{\delta \searrow 0}\int\limits_{O'}\int\limits_{\B(x,r)} \big(\Q_u(x,x')\big)^p \rho_{\delta}(x,x')\D \m(x')\D \m(x)\geq \\
                \geq \lims\limits_{r\searrow 0}\limi\limits_{\delta \searrow 0}\mathcal{I}^L_{E,r}[\rho_{\delta}]\lim\limits_{n\to +\infty}\int\limits_{O''}\Big(\lip\big[u^{r_n}\big](x)\Big)^p\D \m(x).
            \end{gathered}
        \end{equation}
        For each $n\in \mathbb{N}$ we put $v_n\coloneqq u^{r_n}$. Thus, in connection with \cref{def:ScChEnerg}, we shall show that $(v_n)_{n\in \mathbb{N}}$ converges in $\Leb^p_{Loc}(O'')$ to $u$ as $n\to +\infty$.

        According to \eqref{eq:ApproxSeq0} from \cref{lem:ApproxSeq}, it suffices to show that
        \begin{equation}
            \lims\limits_{r\searrow 0} \int\limits_{O''}\fint\limits_{\B(x,r)} \big(\dX_u(x,x')\big)^p\D \m(x')=0.
        \end{equation}
        To prove that, we notice first that
        \begin{equation}
            \lims\limits_{r\searrow 0}\fint\limits_{\B(x,r)} \big(\dX_u(x,x')\big)^p=0
        \end{equation}
        for $\m$-a.e. $x\in O''$, which follows directly from assertion $\rm i)$ in \cref{prop:StrDoubProp}. Then, via Jensen's inequality, we write for $\m$-a.e. $x\in O''$ that
        \begin{equation}
            \fint\limits_{\B(x,r)}\big(\dX_u(x,x')\big)^p \D\m(x')\leq 2^{p-1}\big|u(x)\big|^p+2^{p-1}\fint\limits_{\B(x,r)} \big|u(x')\big|^p \D \m(x').
        \end{equation}
        We claim the right-hand side above converges in $\Leb^1(O'')$ as  $r\searrow 0$. Indeed, the first term is just a fixed function from $\Leb^1(O'')$, while for the second term it holds due to assertion $\rm ii)$ from \cref{prop:StrDoubProp}, so the claim is valid. We are in position to apply \cref{prop:GenLebTh}, which gives us that the necessary convergence.

        What we have is the following estimate:
        \begin{equation}
        \begin{gathered}
            \frac{1}{C}\E_p^{Sc}[u](O'')  \lims\limits_{r\searrow 0}\limi\limits_{\delta \searrow 0}\mathcal{I}^L_{E,r}[\rho_{\delta}]\leq \\
            \leq \lim\limits_{r\searrow 0} \limi\limits_{\delta\searrow 0}\int\limits_{O'}\int\limits_{\B(x,r)} \big(\Q_u(x,x')\big)^p \rho_{\delta}(x,x')\D \m(x')\D \m(x).
            \end{gathered}
        \end{equation}
        Now we just use the inclusion $V\subseteq O''$ to replace $O''$ with $V$ in the left-sided expression above. After this, we maximize the corresponding expression over $\m$-measurable sets $E\subseteq \X$ with $V\Subset E$ and juxtapose the result with \eqref{eq:LowIntMolFamily}. This immediately gives \eqref{eq:LowBoundLem0}, which concludes the proof.
    \end{proof}
\end{lem}

Within the next theorem, the lower-bound part of \cref{theo:Res} is established, both in the localized and global forms. For the proof, we recall \cref{lem:ChEnerMeas}, \cref{prop:LatSupExpr}, and \cref{def:LowAdmFam}.
\begin{theo}\label{theo:LowBoundTh}
    Let $\Omega\subseteq \X$ be an open set, suppose $\X$ is strongly doubling within $\Omega$. Let $p\in \R_{\geq 1}$, let $\rho_{\bullet}$ be a lower-admissible family for $\Omega$. Then there exists a constant $C^L\coloneqq C^L[\Omega]\in \R_{>0}$ depending only on $p$, $C_{\mathrm{SD}}[\Omega]$, $C_{\mathrm{M}}^L[\Omega,\rho_{\bullet}]$ such that the following holds. Let $(\Y,\dY)$ be a metric space, let $O\subseteq \Omega$ be an open set, let $f\in \Me(O,\Y)$. Then, given any open set $O'\Subset O$, one has
        \begin{equation}\label{eq:LowBoundTh0}
            C^L\E_p[f](O')\leq \lim\limits_{r\searrow 0}\limi\limits_{\delta \searrow 0}\int\limits_{O'}\int\limits_{\B(x,r)} \big(\Q_f(x,x')\big)^p\rho_{\delta}(x,x')\D \m(x')\D \m(x).
        \end{equation}
        Consequently, the estimate below takes place:
        \begin{equation}\label{eq:LowBoundTh00}
            C^L\E_p[f](O)\leq \limi\limits_{\delta \searrow 0}\int\limits_{O\times O}\big(\Q_f(x,x')\big)^p\rho_{\delta}(x,x')\D (\m\otimes \m)(x,x').
        \end{equation}
    \begin{proof}
        Let $C'=C'(\Omega)$ denote the constant guaranteed by \cref{lem:ApproxSeq}. And let us then define $C^L[\Omega]\coloneqq \dfrac{C_{\mathrm{M}}^L[\Omega,\rho_{\bullet}]}{C'}$.
        
        Fix an open set $O'\Subset O$. In correspondence with \cref{prop:LatSupExpr}, we fix a pairwise disjoint family $(O_j)_{j\in \mathbb{N}}$ of open subsets of $O'$ together with a family $(\phi_j)_{j\in \mathbb{N}}\subseteq \BLIP_1(\Y)$. Put then $u_j\coloneqq \phi_j\circ f$ for each $j\in \mathbb{N}$. By trivial reasons, we have $O_j\Subset O$ and $u_j\in \Leb^p_{Loc}(O_j)$ for each $j\in \mathbb{N}$.

        Fix $j\in \mathbb{N}$. From \cref{lem:LowBoundLem}, together with \eqref{eq:LowAdmIneq} from \cref{def:LowAdmFam}, we can derive that
        \begin{equation}
            C^L[\Omega]\E^{Sc}_p[u_j](V) \leq \lim\limits_{r\searrow 0}\limi\limits_{\delta \searrow 0}\int\limits_{O_j}\int\limits_{\B(x,r)}\big(\Q_{u_j}(x,x')\big)^p \rho_{\delta}(x,x')\D \m(x')\D \m(x)
        \end{equation}
       for any compact set $V\subseteq O_j$. Taking the supremum from the inequality above over such sets $V$, with \eqref{eq:ChEnerMeas} from \cref{lem:ChEnerMeas} in mind, we get
        \begin{equation}\label{eq:LowBoundTh1}
            C^L[\Omega] \E^{Sc}_p[u_j](O_j)\leq \lim\limits_{r\searrow 0}\limi\limits_{\delta \searrow 0}\int\limits_{O_j}\int\limits_{\B(x,r)}\big(\Q_{f}(x,x')\big)^p \rho_{\delta}(x,x')\D \m(x')\D \m(x),
        \end{equation}
        where we also used the fact that $\phi_j$ is $1$-Lipschitz to replace $u_j$ with $f$ on the right.
        
        We now sum up the inequalities as in \eqref{eq:LowBoundTh1} over $j\in \mathbb{N}$ and apply Fatou's lemma for series two times. Exploiting then \eqref{eq:LatSupExpr} from \cref{prop:LatSupExpr} in combination with the arbitrariness of the families $(O_j)_{j\in \mathbb{N}}$ and $(\phi_j)_{j\in \mathbb{N}}$, we come exactly to the estimate in \eqref{eq:LowBoundTh0}.

 To achieve \eqref{eq:LowBoundTh00}, it remains only to take there the supremum over all open sets $O'\Subset O$, apply simple properties of the measure $\E_p[f]$, and use some obvious inequalities.

        The proof is complete.
    \end{proof}
\end{theo}

\subsection{Upper bound}\label{ss:UpBound}
The goal of this section is to prove the upper-bound estimate in \cref{theo:Res}.

We start by proving a pretty abstract statement that allows one to pass from a certain pointwise inequality on a given map to an estimate for integrals involving mollifiers. For the proof, we recall the notation preceding \cref{def:UpAdmFam}.
\begin{prop}\label{prop:AbstrUpIneq}
    Let $p\in \R_{\geq 1}$, let $\rho\in \Me_+\big(\X\times \X\big)$. Let $E_0\subseteq \X$ be a $\m$-measurable set, let $(\Y,\dY)$ be a metric space, let $f\in \Me(E_0,\Y)$. Let $E\Subset E_0$ be a $\m$-measurable set. Suppose there is $R\in \R_{>0}$ with $\B(E,R)\subseteq E_0$ and there exists a function $G\colon E_0\times (0, R]\to \overline{\R}_{\geq 0}$ satisfying the following properties:
    \begin{enumerate}
        \item[$\rm i)$] for every $r\in (0,R]$, the function $G(\,\cdot\,, r)$ is $\m$-measurable;
        \item[$\rm ii)$] there is $C\in \R_{>0}$ such that, given any $r\in (0,R]$, it holds, for $\m$-a.e. $x\in E_0$, that
        \begin{equation}
            \sup\limits_{r'\in \left[\frac{r}{2},r\right]} G(x,r')\leq C G(x,r);
        \end{equation}
        \item[$\rm iii)$] there is $D\subseteq E_0\times E_0$ with $(\m\otimes \m)\big((E_0\times E_0)\backslash D\big)=0$ such that it holds, for any $(x,x')\in D$ with $x\in E$ and $x'\in \big(\B(x,R)\backslash \{x\}\big)$, that
    \begin{equation}
        \big(\Q_f(x,x')\big)^p\leq G\big(x,\dX(x,x')\big)+G\big(x',\dX(x,x')\big).
    \end{equation}
 \end{enumerate}
    Then it holds, for any $r\in (0,R]$, that
    \begin{equation}\label{eq:AbstrUpIneq0}
    \begin{gathered}
        \int\limits_E \int\limits_{\B(x,r)} \big(\Q_f(x,x')\big)^p\rho(x,x')\D \m(x')\D\m(x)\leq \\
        \leq C \mathcal{I}_{E,r}^U[\rho]\sup\limits_{r'\in (0,r]} \int\limits_{E_0} G(x, r')\D \m(x).
        \end{gathered}
    \end{equation}

    \begin{proof}
We fix any $r\in (0,R]$, after which we pick any $\tau_{\bullet}\in \P^U_r$. Put $t_k\coloneqq \sup(\tau_k)$ for each $k\in \mathbb{N}_0$. From the very definition of $\P_r^U$, we know that $\tau_k\subseteq \big[\tfrac{t_k}{2},t_k]$ for each $k\in \mathbb{N}_0$.

Fix $k\in \mathbb{N}_0$. By combining the inequalities from properties $\rm ii)$ and $\rm iii)$ in the premise, we easily derive that, for $(\m\otimes \m)$-a.e. $(x,x')\in E_0\times E_0$ with $x\in E$ and $x'\in \A(x,\tau_k)$, one has the estimate
\begin{equation}
    \big(\Q_f(x,x')\big)^p\leq C\Big(G(x,t_k)+G(x',t_k)\Big).
\end{equation}
We then integrate this estimate over the corresponding subset of $E_0\times E_0$, which is possible due to property $\rm i)$ in the premise, and carry out some simple manipulations with the resulting integral, which gives us the following:
 \begin{equation}\label{eq:AbstrUpIneq1}
            \begin{gathered}
                \int\limits_E \int\limits_{\A(x,\tau_k)} \big(\Q_f(x,x')\big)^p\rho(x,x')\D \m(x')\D\m(x)\leq\\
                \leq C\int\limits_E \int\limits_{\A(x,\tau_k)} \Big(G(x,t_k)+G(x',t_k)\Big)\rho(x,x')\D \m(x')\D\m(x)\leq\\
                \leq C\int\limits_{\B(E,t_k)} G(x,t_k)\int\limits_{\A(x,\tau_k)} \Big(\rho(x,x')+\rho(x',x)\Big)\D \m(x')\D\m(x)\leq\\
                \leq C \biggg(\esup\limits_{x\in \B(E,t_k)}\int\limits_{\A(x,\tau_k)}\Big(\rho(x,x')+\rho(x',x)\Big)\D \m(x')\biggg) \sup\limits_{r'\in (0,r]}\int\limits_{\B(E,r)} G(x,r')\D \m(x).
            \end{gathered}
        \end{equation}
        
Now we sum up the inequalities as in \eqref{eq:AbstrUpIneq1} over $k\in \mathbb{N}_0$, with taking into account that $\big(\Q_f(x,x)\big)^p\rho(x,x)=0$ for any $x\in E_0$ according to our conventions. Recall also that $\B(E,r)\subseteq E_0$. With all this and with \eqref{eq:UpIntMolPart} in mind, we get
\begin{equation}
\begin{gathered}
\int\limits_E \int\limits_{\B(x,r)} \big(\Q_f(x,x')\big)^p\rho(x,x')\D \m(x')\D\m(x)\leq \\
\leq C\mathcal{I}^U_{E,r}[\rho,\tau_{\bullet}]\sup\limits_{r'\in (0,r']}\int\limits_{E_0} G(x,r')\D \m(x).
    \end{gathered}
\end{equation}
It remains only to minimize the right-sided expression over $\tau_{\bullet}\in \P^U_r$, which, together with \eqref{eq:UpIntMolFul}, leads to the desired estimate from \eqref{eq:AbstrUpIneq0}. This finishes the proof.
    \end{proof}
    
\end{prop}

In the next proposition we provide conditions guaranteeing the presence of a pointwise inequality of the form as in the proposition above. Namely, it is the well-known pointwise inequality that is implied by the doubling and Poincar{\'e} inequalities. Precisely as in \cite[Lemma 3.3]{O25}, we formulate it in terms of the Riesz potentials, which is the most convenient way to deal simultaneously with any $p\in \R_{\geq 1}$. But in contrast with the version from \cite[Lemma 3.3]{O25}, here we need to derive also a slightly enhanced variation of the inequality, which is, in fact, achievable almost straightforward from the weaker one.
\begin{prop}\label{prop:PointIneq}
    Let $\Omega\subseteq \X$ be an open set, let $p\in \R_{\geq 1}$, suppose $\X$ is doubling around $\Omega$ and supports a $p$-Poincar{\'e} inequality around $\Omega$. Then there exists a constant $C=C(\Omega)$ depending only on $p$, $C_{\mathrm{D}}[\Omega]$, $C_{\mathrm{P}}[\Omega]$ and there is a parameter $\lambda=\lambda(\Omega)\in \R_{\geq 1}$ such that the following holds. Let $(\Y,\dY)$ be a metric space, let $O\subseteq \X$ be an open set, let $f\in \Me(O,\Y)$. Then there is $Q_f\subseteq O$ with $\m\big(O\backslash Q_f\big)=0$ and, given any $S\Subset O$ with $S\subseteq \Omega$, there is a radius $R=R(\Omega,O,S)\in \R_{>0}$ with $\B\big(S,2\lambda R\big)\subseteq O$ with the property that it holds, for any $r\in (0,R]$, any $x\in S\cap Q_f$, and any $x'\in \B(x,r)\cap Q_f$, that
    \begin{equation}\label{eq:PointIneq0}
        \big(\dX_f(x,x')\big)^p\leq C r^p\Big(\mathcal{R}_p[f]\big(x,\lambda r\big)+\mathcal{R}_p[f]\big(x',\lambda r\big)\Big).
    \end{equation}
    Consequently, with $S$ and $R$ as above, for any $x\in S\cap Q_f$ and any $x'\in \B(x,R)\cap Q_f$, one has
    \begin{equation}\label{eq:PointIneq00}
        \big(\Q_f(x,x')\big)^p\leq  C\bigg(\mathcal{R}_p[f]\Big(x,\lambda \dX(x,x')\Big)+\mathcal{R}_p[f]\Big(x',\lambda \dX(x,x')\Big)\bigg).
    \end{equation}
    \begin{proof}
        The first part of the statement, namely the inequality in \eqref{eq:PointIneq0}, is essentially established in \cite[Lemma 3.3]{O25}. And even though both the doubling and Poincar{\'e} conditions considered there slightly differ from the currently adopted ones, we claim that the proof of the modified version is basically the same. Thus, we omit details related to this part.

        In order to derive the second part of the statement from the first one, we act as follows. Take any $x\in S\cap Q_f$ and $x'\in \B(x,R)\cap Q_f$. If $x=x'$, then \eqref{eq:PointIneq00} holds trivially, so we assume that $x\neq x'$. Putting $r\coloneqq \dX(x,x')$, we plug $r,x,x'$ into \eqref{eq:PointIneq0}, which gives exactly the inequality in \eqref{eq:PointIneq00}. This finishes the proof.
    \end{proof}
\end{prop}

In the next lemma, we do nothing but combine \cref{prop:AbstrUpIneq} and \cref{prop:PointIneq}, which gives us a general estimate from below on the energy of a given map in terms of an integral with a mollifier. In fact, the resulting statement is useful to compare with \cite[Theorem 1]{BBM01}, as the former one is a natural extension of the latter one to the non-smooth situation. For the proof, we recall \cref{prop:RieszBorel} and \cref{prop:DoubProp}, together with the notation prior to \cref{def:UpAdmFam}.
\begin{lem}\label{lem:UpBoundLem}
    Let $\Omega\subseteq \X$ be an open set, let $p\in \R_{\geq 1}$, suppose $\X$ is doubling around $\Omega$ and supports a $p$-Poincar{\'e} inequality around $\Omega$. Then there exists a constant $C=C(\Omega)\in \R_{>0} $ depending only on $p$, $C_{\mathrm{D}}[\Omega]$, $C_{\mathrm{P}}[\Omega]$ such that the following holds. Let $\rho\in  \Me_+\big(\X\times \X\big)$. Let $(\Y,\dY)$ be a metric space, let $O\subseteq \X$ be an open set, let $f\in \Me(O,\Y)$, let $E\Subset O$ be a $\m$-measurable set with $E\subseteq \Omega$. Then there is a radius $R=R(\Omega,O,E)\in \R_{>0}$ such that it holds, for any $r\in (0,R]$, that
        \begin{equation}\label{eq:UpBoundLem0}
            \int\limits_E \int\limits_{\B(x,r)} \big(\Q_f(x,x')\big)^p\rho(x,x')\D\m(x')\D \m(x)\leq C\mathcal{I}_{E,r}^U[\rho]\E_p[f](O).
        \end{equation}
        Consequently, given any $\rho_{\bullet}\subseteq \Me_+\big(\X\times \X\big)$, one has
         \begin{equation}\label{eq:UpBoundLem00}
            \lim\limits_{r\searrow 0} \lims\limits_{\delta\searrow 0}\int\limits_E \int\limits_{\B(x,r)} \big(\Q_f(x,x')\big)^p\rho_{\delta}(x,x')\D\m(x')\D \m(x)\leq C \mathcal{I}^U_E[\rho_{\bullet}] \E_p[f](O).
        \end{equation}
        \begin{proof}
            First of all, we clearly can assume that $\E_p[f](O)<+\infty$, since otherwise \eqref{eq:UpBoundLem0} is in place by trivial reasons.

            Now, let $C'=C'(\Omega)$ be the constant, let $\lambda'=\lambda'(\Omega)$ be the parameter, and let $R'=R'(\Omega,O,E)$ be the radius guaranteed by \cref{prop:AbstrUpIneq}. Then, we can also find $R_0\in (0,R']$ with $\B(E,2\lambda' R_0)\subseteq O$. Define $R\coloneqq \tfrac{1}{\lambda'}\min\big\{R_0,R_{\mathrm{D}}[\Omega]\big\}$ and put $E_0\coloneqq \B(E,\lambda R)$. With this, we define a function $G\colon E_0\times (0,R]\to \overline{\R}_{\geq 0}$ by
            \begin{equation}
                G(x,t)\coloneqq C'\mathcal{R}_p[f]\big(x,\lambda' t).
            \end{equation}

            Assume for the moment that $G$ satisfies all the properties required in the premise of \cref{prop:AbstrUpIneq}. Let also $C_G$ denote the corresponding constant from property $\rm ii)$ there and assume further that this constant depends only on $p$, $C_{\mathrm{D}}[\Omega]$. With this, let us now complete the proof of the first part of the current lemma. 

Fix any $r\in (0,R]$. Applying \eqref{eq:AbstrUpIneq0} from \cref{prop:AbstrUpIneq} and \eqref{eq:DoubProp2} from assertion $\rm ii)$ in \cref{prop:DoubProp}, we can write
            \begin{equation}
                \begin{gathered}
                    \int\limits_E\int\limits_{\B(x,r)}\big(\Q_f(x,x')\big)^p\rho(x,x')\D\m(x')\D \m(x)\leq \\
                    \leq C' C_G\mathcal{I}^U_{E,r}[\rho]\sup\limits_{r'\in (0,r]}\int\limits_{\B(E,r)} \mathcal{R}_p[f]\big(x,\lambda' r\big)\D \m(x)\leq \\
                    \leq C' C_G C_{\mathrm{D}}[\Omega]\mathcal{I}^U_{E,r}[\rho]\E_p[f](O).
                \end{gathered}
            \end{equation}
            Putting $C\coloneqq C' C_G C_{\mathrm{D}}[\Omega]$, we exactly get \eqref{eq:UpBoundLem0}, as desired.

            We move to verifying the assumptions made about $G$. As for property $\rm i)$, it is true immediately by \cref{prop:RieszBorel}. In turn, property $\rm ii)$, as well the fact that $C_G$ can be chosen so that it depends only on $p$, $C_{\mathrm{D}}[\Omega]$, are nothing but assertion $\rm i)$ in \cref{prop:DoubProp}. Finally, the validity of property $\rm iii)$ follows directly from \cref{prop:PointIneq} and particularly from \eqref{eq:PointIneq00} therein.

             To derive \eqref{eq:UpBoundLem00} one just needs to pass in \eqref{eq:UpBoundLem0} first to the upper limit as $\delta \searrow 0$ and then to the lower limit as $r\searrow 0$, after which to use \eqref{eq:UpIntMolFamily}. Thus, we can conclude the proof.
        \end{proof}
\end{lem}

The role of the next simple proposition is to show that under the integrability assumption on a given map, one can reduce all the relevant reasoning involving mollifiers to the local consideration.
\begin{prop}\label{prop:Localization}
Let $\Omega\subseteq \X$ be an open set, let $p\in \mathbb{R}_{\geq 1}$. Let $\rho_{\bullet}$ be a family satisfying the $p$-decay condition for $\Omega$. Let $(\Y, \dY)$ be a metric space, let $O\subseteq \Omega$ be an open set, let $f\in \Leb^p(O, \Y)$. Then one has the following:
    \begin{equation}\label{eq:Localization0}
\lim\limits_{r\searrow 0}\lims\limits_{\delta\searrow 0} \int\limits_{O}\int\limits_{O\cap \A(x,\infty,r)}\big(\Q_f(x,x')\big)^p\rho_{\delta}(x,x')\D \mathfrak{m}(x')\D \mathfrak{m}(x)=0.
\end{equation}
\begin{proof}
    Essentially, the proof consists of the following steps. Pick $y\in \Y$ with $\dY\big(y,f(\,\cdot\,)\big)\in \Leb^p(O)$. Apply the triangle and Jensen's inequalities to write
    \begin{equation}
        \big(\Q_f(x,x')\big)^p\leq 2^{p-1}\frac{\Big(\dY\big(f(x),y\big)\Big)^p}{\dX(x,x')}+2^{p-1}\frac{\Big(\dY\big(y,f(x')\big)\Big)^p}{\dX(x,x')}
    \end{equation}
    for all distinct $x,x'\in O$. Then, plug this into the integral in \eqref{eq:Localization0} and perform some simple manipulations with the result to separate the function and the mollifiers. It remains only to use \eqref{eq:MolDecCond} from \cref{def:MolDecCond}. The more detailed proof can be found in \cite[Proposition 3.1]{O25}.
\end{proof}
\end{prop}

The theorem below establishes the upper-bound inequality from \cref{theo:Res}, together with its localized analog. For the proof, we recall \cref{def:StrExtDom}, \cref{def:UpAdmFam}, and \cref{def:MolDecCond}.
\begin{theo}\label{theo:UpBoundTh}
    Let $\Omega\subseteq \X$ be an open set, let $p\in \R_{\geq 1}$, suppose $\X$ is doubling around $\Omega$ and supports a $p$-Poincar{\'e} inequality around $\Omega$. Let $\rho_{\bullet}$ be an upper-admissible family for $\Omega$. Then there exists a constant $C^U\coloneqq C^U[\Omega]\in \R_{>0}$ depending only on $p$, $C_{\mathrm{D}}[\Omega]$, $C_{\mathrm{P}}[\Omega]$, $C^U_{\mathrm{M}}[\Omega, \rho_{\bullet}]$ such that the following holds. Let $(\Y,\dY)$ be a metric space, let $O\subseteq \Omega$ be an open set, let $f\in \Me(O,\Y)$. Then, given any $\m$-measurable set $E\Subset O$, one has
    \begin{equation}\label{eq:UpBoundTh0}
           \lim\limits_{r\searrow 0}\lims\limits_{\delta\searrow 0}\int\limits_{E}\int\limits_{\B(x,r)} \big(\Q_f(x,x')\big)^p \rho_{\delta}(x,x') \D(\m\otimes \m)(x,x')\leq C^U\E_p[f](O).
        \end{equation}
        If, in addition, it holds that $O$ has the strong $p$-extension property with respect to $\Y$, that $\rho_{\bullet}$ satisfies the $p$-decay condition for $O$, and that $f\in \Leb^p(O,\Y)$, then the estimate below takes place:
        \begin{equation}\label{eq:UpBoundTh00}
            \lims\limits_{\delta\searrow 0}\int\limits_{O\times O} \big(\Q_f(x,x')\big)^p \rho_{\delta}(x,x') \D(\m\otimes \m)(x,x')\leq C^U\E_p[f](O).
        \end{equation}
        \begin{proof}
        Let $C'=C'(\Omega)$ denote the constant guaranteed by \cref{lem:UpBoundLem}. And let us also define $C^U[\Omega]\coloneqq C'C^U_{\mathrm{M}}[\Omega,\rho_{\bullet}]$.

        In order to prove the first part, we just need to combine \eqref{eq:UpBoundLem00} from \cref{lem:UpBoundLem} with \eqref{eq:UpAdmFam} from \cref{def:UpAdmFam}. It remains only to verify the second part.
        
            Since $f\in \Leb^p(O,\Y)$ and since $O$ has the strong $p$-extension properties with respect to $\Y$ by the assumptions, there exist a metric space $(\Z,\dZ)$ and a $\Z$-extension $F$ of $f$ with
            \begin{equation}\label{eq:UpBoundTh1}
                \E_p[F](O)=\lim\limits_{R\searrow 0}\E_p[F]\big(\B(O,R)\big).
            \end{equation}
            
Now we pick any $R\in \R_{>0}$ and apply \cref{lem:UpBoundLem} to the data  $\B(O,R)$, $F$, $O$. This gives us
\begin{equation}
    \lim\limits_{r\searrow 0}\lims\limits_{\delta \searrow 0} \int\limits_O\int\limits_{\B(x,r)} \big(\Q_F(x,x')\big)^p\rho_{\delta}(x,x')\D \m(x')\D \m(x)\leq C^U[\Omega]\E_p[F]\big(\B(O,R)\big).
\end{equation}
Since $F$ is a $\Z$-extension of $f$, we know that $F\in \Leb^p(O,\Z)$, so we are in position to apply \cref{prop:Localization}. Together with the estimate above, this easily implies the following:
\begin{equation}
    \lims\limits_{\delta \searrow 0} \int\limits_{O\times O} \big(\Q_F(x,x')\big)^p\rho_{\delta}(x,x')\D (\m\otimes \m)(x,x')\leq C^U[\Omega]\E_p[F]\big(\B(O,R)\big).
\end{equation}
            Now we are able to pass in the above inequality to the limit as $R\searrow 0$, after which we use \eqref{eq:UpBoundTh1}. Finally, to conclude the proof, we just notice again that $F$ is a $\Z$-extension of $f$, whence we can replace $F$ with $f$ in both sides of the resulting inequality, which leads exactly to the estimate in \eqref{eq:UpBoundTh00}, as desired. 
        \end{proof}
\end{theo}

\subsection{Simplified assumptions}\label{ss:SimplAssMol}
The purpose of this section is to present adequately verifiable sufficient conditions for a family of mollifiers to be admissible.

\begin{rema}
    The way we formulate the conditions in the statements below is chosen so as to highlight similarities and differences between those and the original BBM ones, despite that this may seem somewhat artificial. A detailed comment on this will be made at the end of the section.
\end{rema}

In the lemma below we provide a convenient replacement for the condition from \cref{def:LowAdmFam}.
\begin{lem}\label{lem:SimplLowAdm}
    Let $\Omega\subseteq \X$ be an open set, suppose $\X$ is strongly doubling within $\Omega$. Let $\rho_{\bullet}\subseteq \Me_+\big(\X\times \X\big)$. Suppose there exists a family $\big(\varrho^-_{\delta}\big)_{\delta\in (0,1)}$ of $\mathcal{L}^1$-measurable functions $\R_{> 0}\to \overline{\R}_{\geq 0}$ satisfying the following properties:
    \begin{enumerate}
        \item[$\rm i)$] one has
        \begin{equation}\label{eq:SimplLowAdm0}
            C^-\coloneqq \lim\limits_{r\searrow 0} \limi\limits_{\delta \searrow 0} \int\limits_{(0,r]} \frac{\varrho^-_{\delta}(t)}{t} \D \mathcal{L}^1(t)>0;
        \end{equation}
        \item[$\rm ii)$] given any compact set $V\subseteq \Omega$, there is a radius $R\in \R_{>0}$ with $\B(V,2R)\subseteq \Omega$ such that it holds, for any $\delta\in (0,1)$, that
        \begin{equation}\label{eq:SimplLowAdm00}
            \rho_{\delta}(x,x')\geq \frac{\varrho^-_{\delta}\big(\dX(x,x')\big)}{\m\Big(\B\big(x,\dX(x,x')\big)\Big)}
        \end{equation}
        for $(\m \otimes \m)$-a.e. $(x,x')\in \X\times \X$ with $x\in \B(V,R)$ and $x'\in \big(\B(x,R)\backslash \{x\}\big)$;
        \item[$\rm iii)$] there is $\wbar{r}\in \R_{>0}$ such that it holds, for any $\delta\in (0,1)$ and any $k\in \mathbb{N}_0$, that $\varrho^-_{\delta}$ is constant on $\big(\tfrac{\wbar{r}}{2^{k+1}},\tfrac{\wbar{r}}{2^k}\big]$.
    \end{enumerate}
    Then $\rho_{\bullet}$ is lower-admissible for $\Omega$ and, moreover, the related constant $C_{\mathrm{M}}^L[\Omega,\rho_{\bullet}]$ can be chosen so that it depends only on $C_{\mathrm{SD}}[\Omega]$, $C^-$.
    \begin{proof}
For every $\delta\in (0,1)$ and for each $k\in \mathbb{N}_0$, let $\sigma^-_{\delta,k}$ denote the value of $\varrho_{\delta}^-$ on $\big(\tfrac{\wbar{r}}{2^{k+1}},\tfrac{\wbar{r}}{2^k}\big]$.

Fix a compact set $V\subseteq \Omega$ and let $R$ denote the corresponding radius from property $\rm ii)$, after which we put $R_0\coloneqq \min\big\{R, R_{\mathrm{SD}}[\Omega,V]\big\}$ and $E\coloneqq \B(V,R_0)$. For each $k\in \mathbb{N}_0$ put $r_k\coloneqq \tfrac{\wbar{r}}{2^k}$. Find $k_0\in \mathbb{N}_0$ with $r_{k_0}\leq R_0$. Fix now any $k\in \mathbb{N}_0$ with $k_0\leq k$ and put $t_{k'}\coloneqq \frac{r_k}{2^{k'}}$ for each $k'\in \mathbb{N}_0$. We put also $\tau_{k'}\coloneqq \big(t_{k'+1},t_{k'}\big]$ for each $k'\in \mathbb{N}_0$. And hence we have $\tau_{\bullet}\in \P^L_{r_k}$.

Fix $\delta\in (0,1)$. By property $\rm ii)$ we can find $Q\subseteq E\cap \supp(\m)$ with $\m\big(E\backslash Q\big)=0$ such that
\begin{equation}
    \varrho_{\delta}^-\big(\dX(x,x')\big)\leq \rho_{\delta}(x,x') \m\Big(\B\big(x,\dX(x,x')\big)\Big)
\end{equation}
for any $x\in Q$ and for $\m$-a.e. $x'\in \big(\B(x,R_0)\backslash \{x\}\big)$. This, together with the strong doubling inequality, allows us to write, for each $k'\in \mathbb{N}_0$ and for any $x\in Q$, that
\begin{equation}
    \m\big(\A(x,\tau_{k'})\big)\einf\limits_{x'\in \A\bg(x,\tau_{k'}\bg)}\rho_{\delta}(x,x')\geq \sigma^-_{\delta,k'} \frac{\m\big(\A(x,\tau_{k'})\big)}{\m\big(\B(x,t_{k'})\big)}\geq  \frac{\sigma^-_{\delta,k'}}{C_{\mathrm{SD}}[\Omega]}.
 \end{equation}
 Summing the inequalities as above over all $k'\in \mathbb{N}_0$, with the use of property $\rm iii)$, we can get
 \begin{equation}
 \begin{gathered}
     \sum\limits_{k'\in \mathbb{N}_0}\einf\limits_{x\in E} \Bigg(\m\big(\A(x,\tau_{k'})\big)\einf\limits_{x'\in \A\bg(x,\tau_{k'}\bg)}\rho_{\delta}(x,x')\Bigg)\geq \\
     \geq \frac{1}{C_{\mathrm{SD}}[\Omega]}\sum\limits_{k'\in \mathbb{N}_0}\sigma^-_{\delta,k'}=\frac{1}{\ln(2)C_{\mathrm{SD}}[\Omega]} \int\limits_{(0,r_k]} \frac{\varrho_{\delta}^-(t)}{t}\D \mathcal{L}^1(t).
     \end{gathered}
 \end{equation}
Now we look at \eqref{eq:LowAdmFam}, particularly at all the optimization procedures made there, and juxtapose it with the estimate above. It is not hard to see then that the constant $C_{\mathrm{M}}^L[\Omega,\rho_{\bullet}]\coloneqq \dfrac{C^-}{\ln(2)C_{\mathrm{SD}}[\Omega]}$ indeed makes the condition from \cref{def:LowAdmFam} valid. The proof is complete.
    \end{proof}
\end{lem}

Now we formulate a more transparent analog of the condition from \cref{def:UpAdmFam}.
\begin{lem}\label{lem:SimplUpAdm}
    Let $\Omega\subseteq \X$ be an open set, suppose $\X$ is doubling around $\Omega$. Let $\rho_{\bullet}\subseteq \Me_+\big(\X\times \X\big)$. Suppose there exists a family $\big(\varrho^+_{\delta}\big)_{\delta\in (0,1)}$ of $\mathcal{L}^1$-measurable functions $\R_{> 0}\to \overline{\R}_{\geq 0}$ satisfying the following properties:
    \begin{enumerate}
        \item[$\rm i)$] one has
        \begin{equation}\label{eq:SimplUpAdm0}
            C^+\coloneqq \lim\limits_{r\searrow 0} \lims\limits_{\delta \searrow 0} \int\limits_{(0,r]} \frac{\varrho^+_{\delta}(t)}{t} \D \mathcal{L}^1(t)<+\infty;
        \end{equation}
        \item[$\rm ii)$] there is $R\in \R_{>0}$ such that it holds, for any $\delta\in (0,1)$, that
        \begin{equation}\label{eq:SimplUpAdm00}
            \rho_{\delta}(x,x')\leq \frac{\varrho^+_{\delta}\big(\dX(x,x')\big)}{\m\Big(\B\big(x,\dX(x,x')\big)\Big)}
        \end{equation}
        for $(\m \otimes \m)$-a.e. $(x,x')\in \X\times \X$ with $x\in \B(\Omega,R)$ and $x'\in \big(\B(x,R)\backslash \{x\}\big)$;
        \item[$\rm iii)$] there is $\wbar{r}\in \R_{>0}$ such that it holds, for any $\delta\in (0,1)$ and any $k\in \mathbb{N}_0$, that $\varrho^+_{\delta}$ is constant on $\big(\tfrac{\wbar{r}}{2^{k+1}},\tfrac{\wbar{r}}{2^k}\big]$.
    \end{enumerate}
    Then $\rho_{\bullet}$ is upper-admissible for $\Omega$ and, moreover, the related constant $C_{\mathrm{M}}^U[\Omega,\rho_{\bullet}]$ can be chosen so that it depends only on $C_{\mathrm{D}}[\Omega]$, $C^+$.
    \begin{proof}
       For every $\delta\in (0,1)$ and for each $k\in \mathbb{N}_0$, let $\sigma^+_{\delta,k}$ denote the value of $\varrho_{\delta}^+$ on $\big(\tfrac{\wbar{r}}{2^{k+1}},\tfrac{\wbar{r}}{2^k}\big]$.

Put $R_0\coloneqq \min\big\{R, R_{\mathrm{D}}[\Omega]\big\}$. For each $k\in \mathbb{N}_0$ we put $r_k\coloneqq \tfrac{\wbar{r}}{2^k}$. Find $k_0\in \mathbb{N}_0$ with $r_{k_0}\leq R_0$. Fix now any $k\in \mathbb{N}_0$ with $k_0\leq k$ and put $t_{k'}\coloneqq \frac{r_k}{2^{k'}}$ for each $k'\in \mathbb{N}_0$. We put also $\tau_{k'}\coloneqq \big(t_{k'+1},t_{k'}\big]$ for each $k'\in \mathbb{N}_0$. It then follows that $\tau_{\bullet}\in \P^U_{r_k}$.

Fix $\delta\in (0,1)$. From property $\rm ii)$ it follows that we can find $Q\subseteq \B(\Omega,R)\cap \supp(\m)$ with $\m\big(\B(\Omega,R)\backslash Q\big)=0$ such that
\begin{equation}
    \varrho_{\delta}^+\big(\dX(x,x')\big)\geq \rho_{\delta}(x,x') \m\Big(\B\big(x,\dX(x,x')\big)\Big)
\end{equation}
for any $x\in Q$ and for $\m$-a.e. $x'\in \big(\B(x,R_0)\backslash \{x\}\big)$. This, together with the doubling inequality, allows us to write, for each $k'\in \mathbb{N}_0$ and for any $x\in Q$, the following:
\begin{equation}
\begin{gathered}
    \int\limits_{\A\bg(x,\tau_{k'}\bg)} \Big(\rho_{\delta}(x,x')+\rho_{\delta}(x',x)\Big)\D \m(x')\leq  \\
    \leq 2\sigma^+_{\delta,k'}\frac{\m\big(\A(x,\tau_{k'})\big)}{\m\Big(\B\big(x,t_{k'+1}\big)\Big)}\leq 2C_{\mathrm{D}}[\Omega]\sigma^+_{\delta,k'}.
    \end{gathered}
\end{equation}
Summing the inequalities as above over all $k'\in \mathbb{N}_0$, with the use of property $\rm iii)$, we can get
\begin{equation}
\begin{gathered}
    \sum\limits_{k'\in \mathbb{N}_0} \esup\limits_{x\in \B(\Omega,R)}\int\limits_{\A\bg(x,\tau_{k'}\bg)} \Big(\rho_{\delta}(x,x')+\rho_{\delta}(x',x)\Big)\D \m(x')\leq\\
    \leq 2C_{\mathrm{D}}[\Omega]\sum\limits_{k'\in \mathbb{N}_0} \sigma_{\delta,k'}^+=\frac{2C_{\mathrm{D}}[\Omega]}{\ln(2)}\int\limits_{(0,r_k]}\frac{\varrho_{\delta}^+(t)}{t}\D \mathcal{L}^1(t).
    \end{gathered}
\end{equation}
Now we compare the above estimate with what is in place in \eqref{eq:UpAdmFam}. Once we perform all the optimization procedures therein, we can deduce that the constant $C_{\mathrm{M}}^U[\Omega,\rho_{\bullet}]\coloneqq \dfrac{2C_{\mathrm{D}}[\Omega]C^+}{\ln(2)}$ makes the condition from \cref{def:LowAdmFam} true. The proof is complete.
    \end{proof}
\end{lem}

Finally, we demonstrate a simple situation when the condition from \cref{def:MolDecCond} is fulfilled.
\begin{prop}\label{prop:SimplDecCond}
    Let $\Omega\subseteq \X$ be an open set. Let $p\in \R_{\geq 1}$, let $\rho_{\bullet}\subseteq \Me_+\big(\X\times \X\big)$. Suppose there is a family $\big(\varrho^+_{\delta}\big)_{\delta\in (0,1)}$ of $\mathcal{L}^1$-measurable functions $\R_{>0}\to \overline{\R}_{\geq 0}$ with the following properties:
    \begin{enumerate}
        \item[$\rm i)$] one has
        \begin{equation}\label{eq:SimplDecCond0}
            \lim\limits_{r\searrow 0}\lims\limits_{\delta \searrow 0}\int\limits_{(r,+\infty)}\frac{\varrho_{\delta}^+(t)}{t^{p+1}}\D \mathcal{L}^1(t)=0;
        \end{equation}
        \item[$\rm ii)$] it holds, for any $\delta\in (0,1)$, that
         \begin{equation}\label{eq:SimplDecCond00}
            \rho_{\delta}(x,x')\leq \frac{\varrho^+_{\delta}\big(\dX(x,x')\big)}{\m\Big(\B\big(x,4\dX(x,x')\big)\Big)}
        \end{equation}
        for $(\m\otimes \m)$-a.e. $(x,x')\in \Omega\times \Omega$ with $x\neq x'$;
        \item[$\rm iii)$] there is $\wbar{r}\in \R_{>0}$ such that it holds, for any $\delta\in (0,1)$ and any $k\in \mathbb{Z}$, that $\varrho_{\delta}^+$ is constant on $\big(2^k \wbar{r}, 2^{k+1} \wbar{r}\big]$.
    \end{enumerate}
    Then $\rho_{\bullet}$ satisfies the $p$-decay condition for $\Omega$.
    \begin{proof}
        Fix $\delta\in (0,1)$. For each $k\in \mathbb{Z}$, put $r_k\coloneqq 2^k \wbar{r}$ and let $\sigma_{\delta,k}^+$ be the value of $\varrho_{\delta}^+$ on $\big(r_k, r_{k+1}\big]$.

        By property $\rm ii)$ we can find $Q\subseteq \Omega$ with $\m\big(\Omega\backslash Q\big)=0$ such that it holds, for any $x\in Q$, that
        \begin{equation}
            \max\Big\{\rho_{\delta}(x,x'),\rho_{\delta}(x',x)\Big\}\leq \frac{\varrho^+_{\delta}\big(\dX(x,x')\big)}{\m\Big(\B\big(x,3\dX(x,x')\big)\Big)}
        \end{equation}
        for $\m$-a.e. $x'\in \Omega$ with $x'\neq x$. Picking $k\in \mathbb{Z}$ and $x\in Q$ and putting $A_k\coloneqq \A\big(x,r_{k'+1}, r_{k'}\big)$ for each $k'\in \mathbb{Z}$, we write 
        \begin{equation}
        \begin{gathered}
            \int\limits_{\A\bg(x,\infty,r_{k}\bg)}\frac{\rho_{\delta}(x,x')+\rho_{\delta}(x',x)}{\big(\dX(x,x')\big)^p}\D \m(x')\leq \\
            \leq 2\sum\limits_{\substack{k'\in \mathbb{Z} \\ k'\geq k}} \frac{\sigma^+_{\delta,k'} \m\Big(\A\big(x,r_{k'+1},r_{k'}\big)\Big)}{(r_{k'})^p\m\Big(\B\big(x,3r_{k'}\big)\Big)}\leq \frac{2^{p+1} }{\ln(2)}\int\limits_{\bg(r_{k},+\infty\bg)} \frac{\varrho_{\delta}^+(t)}{t^{p+1}}\D \mathcal{L}^1(t).
            \end{gathered}
        \end{equation}
        By property $\rm i)$, once we pass in the right-sided expression above to the upper limit as $\delta \searrow 0$ and then to the limit as $k\to -\infty$, we come exactly to \eqref{eq:MolDecCond}, as desired. This finishes the proof.
    \end{proof}
\end{prop} \noindent We note that the number four in \eqref{eq:SimplDecCond00} cannot be removed in general. This is essentially due to that we operate in the proof not only on small, but also large scales.

\begin{rema}
    Now we wish to give some clarifying comments on the conditions from the statements of the current section.

    The additional families appearing in \cref{lem:SimplLowAdm} and in \cref{lem:SimplUpAdm} with \cref{prop:SimplDecCond} should be interpreted naturally as certain radial, piecewise constant bounds from below and from above, respectively, for the given family of mollifiers, the fact of which can be extracted from the inequalities in \eqref{eq:SimplLowAdm00}, \eqref{eq:SimplUpAdm00}, and \eqref{eq:SimplDecCond00}. The only correction are the weights added therein, which are necessary, broadly speaking, to suppress the dependence on a ``reference'' point $x$. 
    
    As a continuation to what is said in the remark under \cref{def:UpAdmFam}, we also draw the reader's attention to the obvious similarities between the integrals in \eqref{eq:SimplLowAdm0} and in \eqref{eq:SimplUpAdm0} and the integral in \eqref{eq:OrigNormalCond}. Up to the above-mentioned weights, which bring obvious modifications in the Euclidean case, the former quantities should be noticed to possess basically the same form. A analogous situation is with the constraints given by \eqref{eq:SimplDecCond0} and \eqref{eq:OrigLocalCond}, but with the former one being even a weaker version of the latter one. These circumstances make the corresponding conditions into perfect extensions of the original ones to the non-smooth setting.
    
    The relevant difference between the BBM assumptions given in \eqref{eq:OrigNormalCond} and \eqref{eq:OrigLocalCond} and the ones listed in \cref{lem:SimplLowAdm}, \cref{lem:SimplUpAdm}, and \cref{prop:SimplDecCond} lies in the presence of properties $\rm iii)$ in each of these statements, which impose some special restrictions on the bounding families. Namely, the families there are required to consist of step functions, with the steps, moreover, shrinking at a certain exponential rate, common for all the functions, when approaching zero. This fixed rate could be easily chosen to be different, but this, however, would also change some details of our proofs, as well as some of the related universal constants.
    
    As the final comment here, we note that the presence of the described fixed rate, while not being truly important for the upper-bound part, is vital for the lower-bound part, even if dealt with only one given family. Indeed, if such a rate was not in place, there would not be possible to prove a suitably modified version of \cref{lem:ApproxSeq}, as the universal constant therein would become infinite due to the necessity to use the doubling inequality an unbounded number of times. And this actually gives the main constraint on what mollifiers are covered by our results. Roughly speaking, the essence of these limitations is that there must be a certain control in place on how mollifiers oscillate. Thus, it is still an open question of how close one is able to approach the original BBM conditions, depending on the given assumptions on the underlying space, within the considered problem.
\end{rema}

\subsection{Admissible mollifiers}\label{ss:ExFamMol}
In this section we present some particular examples of mollifiers that fall into the scope of our result. 

There are numerous specific mollifiers of the considered type that typically appear in the literature. Many of them, within the Euclidean setting, are discussed in \cite{BBM01} and \cite{B02}. The corresponding analogs in the singular context can be found in \cite{LPZ22} and \cite{O25}. And here we want basically to duplicate all the families of mollifiers listed therein, but with adding an extra one, namely the counterpart of the family given by \eqref{eq:NonMonFamEx}.

Let $p\in \mathbb{R}_{\geq 1}$. For every $\delta\in (0,1)$, define functions $\rho^{\kappa,p}_{\delta}\in \Me_+\big(\X\times \X\big)$, $\kappa\in \{1,2,3,4,5\}$, by
\begin{gather}
\rho^{1,p}_{\delta}(x,x')\coloneqq\frac{\delta \big(\dX(x,x')\big)^{p\delta} \chi_{\mathrm{B}(x,1)}(x')}{\mathfrak{m}\Big(\mathrm{B}\big(x,4\mathsf{d}(x,x')\big)\Big)},\label{eq:ExMolFam1}\\
\rho^{2,p}_{\delta}(x,x')\coloneqq\frac{\big(\dX(x,x')\big)^p\chi_{\mathrm{B}(x,\delta)}(x')}{{\delta}^p \mathfrak{m}\big(\mathrm{B}(x,\delta)\big)},\label{eq:ExMolFam2}\\
\rho^{3,p}_{\delta}(x,x')\coloneqq\frac{\chi_{\mathrm{B}(x,\delta)}(x')}{\mathfrak{m}\big(\mathrm{B}(x,\delta)\big)},\label{eq:ExMolFam3}\\
\rho^{4,p}_{\delta}(x,x')\coloneqq\frac{\big(\dX(x,x')\big)^p\chi_{\mathrm{B}(x,\delta)}(x')}{{\delta}^p \mathfrak{m}\Big(\mathrm{B}\big(x,\mathsf{d}(x,x')\big)\Big)}\label{eq:ExMolFam4},\\
\rho^{5,p}_{\delta}(x,x')\coloneqq \frac{\chi_{\mathrm{A}(x,1,\delta)}(x')}{\big|\ln(\delta)\big| \mathfrak{m}\Big(\mathrm{B}\big(x,4\mathsf{d}(x,x')\big)\Big)}\label{eq:ExMolFam5},
\end{gather}
respectively. We note that some of the families above possess an independent meaning, such as the one given by \eqref{eq:ExMolFam1}, which is known to be related to the fractional Sobolev seminorms.

\begin{rema}
    The presence of the number four in \eqref{eq:ExMolFam1} and \eqref{eq:ExMolFam5} is due to that without it, the corresponding families may not fulfill the decay condition from \cref{def:MolDecCond}, unless the underlying space satisfies a certain doubling condition, stronger than the one we want to deal with. And in general, one can put there any number greater than two. All this can be extracted from \cref{prop:SimplDecCond} and \cref{lem:MolAdmEx}.
\end{rema}

Now we demonstrate the necessary assumptions under which the mollifiers listed above satisfy the conditions from \cref{ss:ExFamMol}.
\begin{lem}\label{lem:MolAdmEx}
Let $\Omega\subseteq \X$ be an open set, let $p\in \R_{\geq 1}$. Then, for each $\kappa\in \{1,2,3,4,5\}$, the following assertions hold.
\begin{enumerate}
    \item[$\rm i)$] Suppose $\X$ is strongly doubling within $\Omega$. Then $\rho^{\kappa,p}_{\bullet}$ is lower-admissible for $\Omega$ and, moreover, the related constant $C_{\mathrm{M}}^L\big[\Omega, \rho^{\kappa,p}_{\bullet}\big]$ can be chosen so that it depends only on $p$, $C_{\mathrm{SD}}[\Omega]$.
    \item[$\rm ii)$] Suppose $\X$ is doubling around $\Omega$. Then $\rho^{\kappa,p}_{\bullet}$ is upper-admissible and satisfies the $p$-decay condition, both for $\Omega$, and, moreover, the related constant $C_{\mathrm{M}}^U\big[\Omega,\rho^{\kappa,p}_{\bullet}\big]$ can be chosen so that it depends only on $p$, $C_{\mathrm{D}}[\Omega]$.
\end{enumerate}
\begin{proof}
Proofs for all the families are quite similar. Namely, one just needs to find the corresponding families from \cref{lem:SimplLowAdm}, \cref{lem:SimplUpAdm}, and \cref{prop:SimplDecCond}, which is essentially straightforward. The detailed proof, with some minor differences from the current formulation, for the family $\rho^{1,p}_{\bullet}$ is given within \cite[Theorem 4.4]{O25}. Since the main novelty of our manuscript is the validity of \cref{theo:Res} for the family $\rho^{5,p}_{\bullet}$, we provide the full proof only for it. So, for every $\delta\in (0,1)$ we put $\rho_{\delta}\coloneqq \rho^{5,p}_{\delta}$. Put also $C_{\mathrm{SD}}\coloneqq C_{\mathrm{SD}}[\Omega]$ for brevity.

We put $r_k\coloneqq \tfrac{1}{2^k}$ for each $k\in \mathbb{N}_0$ and fix $\delta\in (0,1)$. Let $K_{\delta}$ denote the smallest number $k\in \mathbb{N}_0$ such that $r_{k+1}\not\in (\delta,1]$, put also $r^-_{\delta}\coloneqq r_{K_{\delta}}$ and $r^+_{\delta}\coloneqq r_{K_{\delta}+1}$. It is clear then that
\begin{equation}
    \big|\log_2(2\delta)\big|\leq K_{\delta}\leq \big|\log_2(\delta)\big|.
\end{equation}

We now prove assertion $\rm i)$. For every $\delta\in (0,1)$, define a function $\varrho^-_{\delta}\colon \R_{>0}\to \R_{\geq 0}$ by
\begin{gather}
    \varrho^-_{\delta}(t)\coloneqq \frac{1}{(C_{\mathrm{SD}})^2\big|\ln(\delta)\big|}\chi_{\bg(r^-_{\delta},1\bg]}(t),
\end{gather}
which is clearly $\mathcal{L}^1$-measurable. The family $\big(\varrho_{\delta}^-\big)_{\delta\in (0,1)}$ can be seen to satisfy properties $\rm ii)$ and $\rm iii)$ from \cref{lem:SimplLowAdm}, with the former due to the strong doubling condition. Thus, it holds, for each $k\in \mathbb{N}_0$ and for any $\delta\in (0,1)$ with $K_{\delta}\geq k$, that
\begin{gather}
    \frac{1}{\ln(2)}\int\limits_{(0,r_k]} \frac{\varrho_{\delta}^-(t)}{t}\D \mathcal{L}^1(t)= \frac{K_{\delta}-k}{(C_{\mathrm{SD}})^2\big|\ln(\delta)\big|}\geq \frac{\big|\log_2(2\delta)\big|-k}{(C_{\mathrm{SD}})^2\big|\ln(\delta)\big|}.
\end{gather}
After passing above to the lower limit as $\delta\searrow 0$ and to the limit as $k\to +\infty$, we deduce that property $\rm i)$ from \cref{lem:SimplLowAdm} is satisfied with $C^-\coloneqq \dfrac{1}{(C_{\mathrm{SD}})^2}$. We now get the necessary conclusion by \cref{lem:SimplLowAdm}.

It remains to prove assertion $\rm ii)$. For every $\delta\in (0,1)$, define a function $\varrho^+_{\delta}\colon \R_{>0}\to \R_{\geq 0}$ by
\begin{gather}
    \varrho^+_{\delta}(t)\coloneqq \frac{1}{\big|\ln(\delta)\big|}\chi_{\bg(r^+_{\delta},1\bg]}(t),
\end{gather}
which is again clearly $\mathcal{L}^1$-measurable. The family $\big(\varrho_{\delta}^+\big)_{\delta\in (0,1)}$ is easily seen to satisfy properties $\rm ii)$ and $\rm iii)$ from both \cref{lem:SimplUpAdm} and \cref{prop:SimplDecCond}. Then, it holds, for each $k\in \mathbb{N}_0$ and for any $\delta\in (0,1)$ with $K_{\delta}\geq k$, that
\begin{gather}
    \frac{1}{\ln(2)}\int\limits_{(0,r_k]} \frac{\varrho_{\delta}^+(t)}{t}\D \mathcal{L}^1(t)=\frac{K_{\delta}+1-k}{\big|\ln(\delta)\big|}\leq \frac{\big|\log_2(\delta)\big|+1-k}{\big|\ln(\delta)\big|},\\
    \int\limits_{\bg(r_k,+\infty\bg)}\frac{\varrho_{\delta}^+(t)}{t^{p+1}}\D \mathcal{L}^1(t)\leq \frac{1}{p\big|\ln(\delta)\big|}\Big(\tfrac{1}{r_k}-1\Big).
\end{gather}
After passing in each of the above inequalities to the upper limit as $\delta\searrow 0$ and to the limit as $k\to +\infty$, we infer that property $\rm i)$ from \cref{lem:SimplUpAdm} is satisfied with $C^+\coloneqq 1$ and that property $\rm i)$ from \cref{prop:SimplDecCond} holds. The conclusion thus follows from \cref{lem:SimplUpAdm} and \cref{prop:SimplDecCond}.

The proof is complete.
\end{proof}
\end{lem}

\begin{rema}
    As one can easily observe from \cref{lem:MolAdmEx}, all the mollifiers given via \eqref{eq:ExMolFam1}-\eqref{eq:ExMolFam5} are fully applicable not only to \cref{theo:Res}, but also to the separate statements provided in \cref{ss:LowBound} and \cref{ss:UpBound}.
\end{rema}

As the conclusive statement, we present the following corollary, obtained by combining \cref{def:PoinSpace}, \cref{prop:DoubStrDoub}, \cref{def:AdmFam}, \cref{theo:LowBoundTh}, \cref{theo:UpBoundTh}, and \cref{lem:MolAdmEx}, where we once again formulate the resulting estimates, but for the specific mollifiers listed in this section.
\begin{cor}
    Let $p\in \R_{\geq 1}$, suppose $\X$ is a $p$-Poincar{\'e} space. Then there exists a constant $C\in \R_{>0}$ depending only on $p$, $C_{\mathrm{D}}[\X]$, $C_{\mathrm{P}}[\X]$ such that the following holds. Let $(\Y,\dY)$ be a metric space, let $O\subseteq \X$ be an open set having the strong $p$-extension property with respect to $\Y$, let $f\in \Leb^p(O,\Y)$. Then the estimates below take place for every $\kappa\in \{1,2,3,4,5\}$:
    \begin{equation}
    \begin{gathered}
\frac{1}{C}\mathrm{E}_p[f](O)\leq \limi\limits_{\delta\searrow 0} \int\limits_{O\times O} \big(\Q_f(x,x')\big)^p\rho^{p,\kappa}_{\delta}(x,x') \D (\mathfrak{m}\otimes \mathfrak{m})(x,x') 
\leq  \\
\leq\lims\limits_{\delta\searrow 0} \int\limits_{O\times O} \big(\Q_f(x,x')\big)^p\rho^{p,\kappa}_{\delta}(x,x') \D (\mathfrak{m}\otimes \mathfrak{m})(x,x')  \leq  C\mathrm{E}_p[f](O).
\end{gathered}
\end{equation}
\end{cor}

{\bf Acknowledgments.} The work continues a research started by the author as part of his Bachelor Thesis at the Moscow Institute of Physics and Technology. The study was supported by the RSF under grant №24-11-00170.

\printbibliography
\end{document}